\documentclass[12pt]{article}

\usepackage[dvips]{graphicx}
\usepackage{makeidx}
\usepackage{amsmath}
\usepackage{amsfonts}

\newcommand {\etavec}{\boldsymbol{\eta}}

\newcommand {\sign}{{\mathrm{sign}}}

\newcommand {\xivec}{{\boldsymbol{\xi}}}
\newcommand {\lambdavec}{{\boldsymbol{\lambda}}}
\newcommand {\varevec}{{\boldsymbol{\varepsilon}}}

\newcommand {\thetavec}{{\boldsymbol{\theta}}}

\newfont{\pseudocode}{cmtt10}

\newtheorem{thm}{Theorem}
\newtheorem{cor}{Corollary}
\newtheorem{lem}{Lemma}

\newtheorem{defn}{Definition}
\newtheorem{cond}{Condition}
\newtheorem{assume}{Assumption}

\newcommand{\bfx}{\mathbf{x}}
\newcommand{\bff}{\mathbf{f}}

\newcommand{\bfz}{\mathbf{z}}

\newcommand{\bfh}{\mathbf{h}}

\newcommand{\bfu}{\mathbf{u}}

\newcommand{\bfq}{\mathbf{q}}

\newcommand{\dist}{\mathrm{dist}}

\newcommand{\pd}{{\partial}}

\newcommand{\Real}{\mathbb{R}}
\newcommand{\Natural}{\mathbb{N}}
\newcommand{\Numbers}{\mathbb{Z}}

\newcommand{\norm}[1]{\left\Vert#1\right\Vert}

\newcommand{\norms}[1]{\norm{#1}_{\mathcal{A}}}
\newcommand{\normss}[2]{\norm{#1}_{\mathcal{A}_{#2}}}
\newcommand{\normsinf}[2]{\norm{#1}_{\mathcal{A}_\infty,#2}}

\makeindex

\begin{document}



\title{\bf Non-uniform Small-gain Theorems for Systems with Unstable Invariant Sets}

\author{Ivan Tyukin\thanks{{\bf Corresponding author}. Laboratory for Perceptual Dynamics, RIKEN (Institute for Physical and Chemical Research)
                           Brain Science Institute, 2-1, Hirosawa, Wako-shi, Saitama, 351-0198, Japan, e-mail:
                           tyukinivan@brain.riken.jp}, Erik
                           Steur\thanks{Dept. of Mechanical Engineering, Dynamics and Control, Eindhoven University of Technology, P.O. Box 513
                           5600 MB,  Eindhoven, The Netherlands, e-mail: e.steur@student.tue.nl}, Henk
                           Nijmeijer \thanks{Dept. of Mechanical Engineering, Dynamics and Control, Eindhoven University of Technology, P.O. Box 513
5600 MB,  Eindhoven, The Netherlands, e-mail: h.nijmeijer@tue.nl},
Cees van Leeuwen\thanks{Laboratory for Perceptual Dynamics, RIKEN
(Institute for Physical and Chemical Research) Brain Science
Institute, 2-1, Hirosawa, Wako-shi, Saitama, 351-0198, Japan,
e-mail: ceesvl@brain.riken.jp} }
\date{\small \it Submitted to SIAM Journal of Control and Optimization (10 October 2006)}
\maketitle{}

\begin{abstract}

We consider the problem of asymptotic convergence to invariant
sets in interconnected nonlinear dynamic systems. Standard
approaches often require that the invariant sets be uniformly
attracting. e.g. stable in the Lyapunov sense. This, however, is
neither a necessary requirement, nor is it always useful. Systems
may, for instance, be inherently unstable (e.g. intermittent,
itinerant, meta-stable) or the problem statement may include
requirements that cannot be satisfied with stable solutions. This
is often the case in general optimization problems and in
nonlinear parameter identification or adaptation. Conventional
techniques for these cases rely either on detailed knowledge of
the system's vector-fields or
 require boundeness of its states. The presently proposed method
relies only on estimates of the input-output maps and steady-state
characteristics. The method requires the possibility of
representing the system as an interconnection of a stable,
contracting, and an unstable, exploratory part. We illustrate with
examples how the method can be applied to problems of analyzing
the asymptotic behavior of locally unstable systems as well as to
problems of parameter identification and adaptation in the
presence of nonlinear parametrizations. The relation of our
results to conventional small-gain theorems is  discussed.
\end{abstract}


{ {\bf Keywords:} non-uniform convergence, weakly attracting sets,
small-gain theorems, input-output stability}


\section{Notation}

Throughout the paper we use the following notational conventions.
Symbol $\Real$ denotes the field of real numbers, symbol $\Real_+$
stands for the following subset of $\Real$: $\Real_+=\{x\in\Real|
 \  x\geq 0\}$; $\Natural$ and $\Numbers$ denote the set of natural
numbers and its extension to the negative domain respectively.

Let ${\Omega}$ be a set, by symbol $\mathcal{S}\{\Omega\}$ we
denote the set of all subsets of $\Omega$. Symbol $\mathcal{C}^k$
denotes the space of functions that are  at least $k$ times
differentiable; $\mathcal{K}$ denotes the class of all strictly
increasing functions $\kappa: \Real_+\rightarrow \Real_+$ such
that $\kappa(0)=0$. If, in addition,
$\lim_{s\rightarrow\infty}\kappa(s)=\infty$ we say that
$\kappa\in\mathcal{K}_\infty$. Further, $\mathcal{K}_e$ (or
$\mathcal{K}_{e,\infty}$) denotes the class of functions of which
the restriction to the interval $[0,\infty)$ belongs to
$\mathcal{K}$ (or $\mathcal{K}_\infty$). Symbol $\mathcal{KL}$
denotes the class of functions
$\beta:\Real_+\times\Real_+\rightarrow\Real_+$ such that
$\beta(\cdot,0)\in\mathcal{K}$ and $\beta(0,\cdot)$ is
monotonically decreasing.

Let $\bfx\in\Real^n$ and $\bfx$ can be partitioned into two
vectors $\bfx_1\in\Real^q$, $\bfx_1=(x_{11},\dots,x_{1q})^T$,
$\bfx_2\in\Real^p$, $\bfx_2=(x_{21},\dots,x_{2p})^T$ with $q+p=n$,
then $\oplus$ denotes their concatenation:
$\bfx=\bfx_1\oplus\bfx_2$.

The symbol $\|\bfx\|$ denotes the Euclidian norm in
$\bfx\in\Real^n$. By ${L}^n_\infty[t_0,T]$ we denote the space of
all functions $\bff:\Real_+\rightarrow\Real^n$ such that
$\|\bff\|_{\infty,[t_0,T]}=\sup\{\|\bff(t)\|,t \in
[t_0,T]\}<\infty$, and $\|\bff\|_{\infty,[t_0,T]}$ stands for the
${L}^n_\infty[t_0,T]$ norm of $\bff(t)$. Let $\mathcal{A}$ be a
set in ${\Real^n}$ and $\|\cdot\|$ be the usual Euclidean norm in
$\Real^n$. By the symbol $\norms{\cdot}$ we denote the following
induced norm:
\[
\norms{\bfx}=\inf_{\bfq\in\mathcal{A}}\{\|\bfx-\bfq\|\}
\]
Let $\Delta\in\Real_+$ then the notation $\normss{\bfx}{\Delta}$
stands for the following equality:
\[
\normss{\bfx}{\Delta}=\left\{\begin{array}{ll}\norms{\bfx}-\Delta,
& \norms{\bfx}>\Delta\\
0, & \norms{\bfx}\leq \Delta\end{array}\right.
\]
The symbol $\normsinf{\cdot}{[t_0,t]}$ is defined as follows:
\[
\normsinf{\bfx(\tau)}{[t_0,t]}= \sup_{\tau\in[t_0,t]}
\norms{\bfx(\tau)}
\]

\section{Introduction}

In many fields of science, from systems and control theory to
physics, chemistry, and biology, it is of fundamental importance
to analyze the asymptotic behavior of dynamical systems. Most of
these analyses are based around the concept of Lyapunov stability
\cite{Lyapunov}, \cite{Yoshizawa:1960}, \cite{Vorotnikov}, i.e.
continuity of the flow
$\bfx(t,\bfx_0):\Real_+\times\Real^n\rightarrow
L_\infty^n[t_0,\infty]$ with respect to $\bfx_0$ \cite{Fradkov99},
 in combination with the standard notion of  an
{\it attracting set} \cite{Guckenheimer_2002}:
\begin{defn}\label{defn:attracting_set} A set $\mathcal{A}$ is an attracting set iff it is

i) closed, invariant, and

ii) for some neighborhood $\mathcal{V}$ of $\mathcal{A}$ and for
all $\bfx_0\in\mathcal{V}$ the following conditions hold:
\begin{equation}\label{eq:attracting_set}
\bfx(t,\bfx_0)\in\mathcal{V} \ \forall \ t\geq 0;
\end{equation}
\begin{equation}\label{eq:attracting_limit}
\lim_{t\rightarrow\infty}\norms{\bfx(t,\bfx_0)}=0
\end{equation}
\end{defn}
Condition (\ref{eq:attracting_set}) in Definition
\ref{defn:attracting_set} stipulates the existence of a trapping
region $\mathcal{V}$ which is  a neighborhood of $\mathcal{A}$.
Condition (\ref{eq:attracting_limit}) ensures attractivity, or
convergence to $\mathcal{S}$. Due to condition
(\ref{eq:attracting_set}), convergence to $\mathcal{A}$ is uniform
with respect to $\bfx_0$ in the neighborhood of $\mathcal{A}$,
i.e. every trajectory which starts in $\mathcal{V}$ remains in
$\mathcal{V}$ for $t\geq 0$ and converges to $\mathcal{A}$ at
$t\rightarrow\infty$.

Although the conventional concepts of attracting set and Lyapunov
stability are a powerful tandem in various applications, some
problems cannot be solved within this framework. Condition
(\ref{eq:attracting_set}), for example, could be violated in
systems with intermittent, itinerant, or meta-stable dynamics. In
general the condition does not hold when the system dynamics,
loosely speaking, is exploring rather than contracting. Such
systems appear naturally in the context of global optimization.
For instance, in \cite{Shang_96} finding the global minimum of a
differentiable cost function $Q: \Real^n\rightarrow \Real_+$ in a
bounded subset $\Omega_x\subset\Real^n$ is achieved by splitting
the search procedure into a locally attracting gradient
$\mathcal{S}_a$, and a wandering part  $\mathcal{S}_w$:
\begin{equation}\label{eq:global_search}
\begin{split}
\mathcal{S}_a: &\  \dot{\bfx} =-\mu_x \frac{\pd Q(\bfx)}{\pd \bfx}
+
\mu_t T(t), \ \mu_x, \mu_t\in \Real_+ \\
\mathcal{S}_w: & \ T(t)=h\{t,\bfx(t)\}, \ h:\Real_+\times
L^n_{\infty}[t_0,t]\rightarrow L^n_\infty[t_0,t]
\end{split}
\end{equation}
The trace function, $T(t)$, in (\ref{eq:global_search}) is
supposed to cover (i.e. be dense in) the whole searching domain
$\Omega_x$. Even though the results in \cite{Shang_96} are purely
simulation studies, they illustrate the superior performance of
algorithms (\ref{eq:global_search}) in a variety of benchmark
problems compared to standard local minimizers and classical
methods of global optimization. Abandoning Lyapunov stability is
likewise advantageous in problems of identification and adaptation
in the presence of general nonlinear parameterization
\cite{Tyukin:IFAC_CONGRESS_2005:2}, in manoeuvring and path
searching \cite{Neural_Computation:Suimetsu:2004}, and in decision
making in intelligent systems \cite{van_Leeuwen_and_Raffone},
\cite{van_Leeuwen_Verver_and_Brinkers}. Systems with attracting,
yet unstable invariant sets are relevant for modelling complex
behavior in biological and physical systems \cite{Ashwin:2005}.
Last but not least, Lyapunov-unstable attracting sets are relevant
in problems of synchronization \cite{Bischi:1998},
\cite{Ott:1994}, \cite{Timme_2002}\footnote{See also
\cite{Pogromsky:2003} where the striking difference between stable
and "almost stable" synchronization in terms of the coupling
strengthes for a pair of the Lorenz oscillators is demonstrated
analytically.}.

Even when it is appropriate to consider a system as stable, we may
be limited in our success in meeting the requirement to identify a
proper Lyuapunov function.%
%
This is the case, for instance, when the system's dynamics is only
partially known. Trading stability requirements for the sake of
convergence might be a possible remedy. Known results in this
direction can be found in \cite{Ilchman_97},
\cite{Pomet92}\footnote{In the Examples section, we demonstrate
how explorative dynamics can solve the problem of simultaneous
state and parameter observation for a system which cannot be
transformed into a canonical adaptive observer form
\cite{Bastin88}.}.

In all the cases that are problematic under condition
(\ref{eq:attracting_set}) of Definition \ref{defn:attracting_set},
condition (\ref{eq:attracting_limit}) --
%
convergence of $\bfx(t,\bfx_0)$ to an invariant set $\mathcal{A}$,
is still a requirement that has to be met.
%
In order to treat these cases analytically we shall, first of all,
move from the standard concept of attracting sets in Definition
\ref{defn:attracting_set} to one that does not assume that the
basin of attraction is necessarily a neighborhood of the invariant
set $\mathcal{A}$. In other words we shall allow convergence which
is not uniform in initial conditions.
This requirement is captured by the concept of weak, or Milnor
attraction \cite{Milnor_1985}:
\begin{defn}\label{defn:Milnor_attracting_set}
A set $\mathcal{A}$ is {weakly attracting}, or Milnor attracting
set iff

i) it is closed, invariant and

ii) for some set $\mathcal{V}$ (not necessarily a neighborhood of
$\mathcal{A}$) with {\it strictly positive measure} and  for all
$\bfx_0\in\mathcal{V}$ limiting relation
(\ref{eq:attracting_limit}) holds
\end{defn}

Conventional methods such as La Salle's invariance principle
\cite{LaSalle:1976} or  center manifold theory \cite{Carr:1981}
can, in principle, address the issue of convergence to weak
equilibria. They do so, however, at the expense of requiring
detailed knowledge of the vector-fields of the ordinary
differential equations of the model. When such information is not
available the system can be thought of as a mere interconnection
of input-output maps. Small-gain theorems
\cite{Zames:66},\cite{Jiang_1994} are usually efficient in this
case. These results, however, apply only under the assumption of
stability of each component in the interconnection.

In the present study we aim to find a proper balance between the
generality of input-output approaches \cite{Zames:66},
\cite{Jiang_1994} in the analysis of convergence and the
specificity of the fundamental notions of limit sets and
invariance  that play a central role in \cite{LaSalle:1976},
\cite{Carr:1981}. The object of our study is a class of systems
that can be decomposed into an attracting, or stable, component
$\mathcal{S}_a$ and an exploratory, generally unstable, part
$\mathcal{S}_w$. Typical systems of this class are nonlinear
systems in cascaded form
\begin{equation}\label{eq:cascade_intro}
\begin{split}
\mathcal{S}_a: \ \dot{\bfx}&=\bff(\bfx,\bfz),\\
\mathcal{S}_w: \ \dot{\bfz}&=\bfq(\bfz,\bfx)
\end{split}
\end{equation}
where the zero solution of the $\bfx$-subsystem is asymptotically
stable in the absence of input $\bfz$, and  the state of the
$\bfz$-subsystem are functions of
$\int_{t_0}^t\|\bfx(\tau)\|d\tau$. Even when both subsystems in
(\ref{eq:cascade_intro}) are stable and the $\bfx$-subsystem does
not depend on state $\bfz$, the cascade can still be unstable
\cite{SIAM:Arcak:2002}. We show, however, that for unstable
interconnections (\ref{eq:cascade_intro}), under certain
conditions that involve only input-to-state properties of
$\mathcal{S}_a$ and $\mathcal{S}_w$, there is a set $\mathcal{V}$
in the system state space, such that trajectories starting in
$\mathcal{V}$ remain bounded. The result is formally stated in
Theorem \ref{theorem:non_uniform_small_gain}. In case an
additional measure of invariance is defined  for $\mathcal{S}_a$
(in our case a steady-state characteristic), a weak, Milnor
attracting set emerges. Its location is completely determined by
the zeros of the steady-state response of system $\mathcal{S}_a$.

We demonstrate how this basic result can be used in problems of
design and analysis of control systems and
identification/adaptation algorithms. In particular, we present an
adaptive observer of state and parameter values for  uncertain
systems which cannot be transformed into a canonic adaptive
observer form \cite{Bastin88}. In the Examples section we present
an application of this result to the problem of reconstructing a
dynamic model of neuronal cell activity.

The paper is organized as follows. In Section 3 we formally state
the problem and provide specific assumptions for the class of
systems under consideration. Section 4 contains the main results
of our present study. In Section 5 we provide several corollaries
of the main result that apply to specific problems. Section 6
contains examples, and Section 7 concludes the paper.

\section{Problem Formulation}

Consider a system that can be decomposed into two interconnected
subsystems, $\mathcal{S}_a$ and $\mathcal{S}_w$:
\begin{equation}\label{eq:attracting_0}
\begin{split}
\mathcal{S}_a: \ & (u_a,\bfx_0)\mapsto \bfx(t)\\
\mathcal{S}_w: \ & (u_w,\bfz_0)\mapsto \bfz(t)
\end{split}
\end{equation}
where $u_a\in\mathcal{U}_a\subseteq L_{\infty}[t_0,\infty]$,
$u_w\in\mathcal{U}_w\subseteq L_\infty[t_0,\infty]$ are the spaces
of inputs to $\mathcal{S}_a$ and $\mathcal{S}_w$, respectively
$\bfx_0\in\Real^n$, $\bfz_0\in\Real^m$ represent initial
conditions, and $\bfx(t)\in\mathcal{X}\subseteq
L_\infty^n[t_0,\infty]$, $\bfz(t)\in\mathcal{Z}\subseteq
L_\infty^m[t_0,\infty]$ are the system states.

System $\mathcal{S}_a$ represents the contracting dynamics. More
precisely, we require that $\mathcal{S}_a$ is input-to-state
stable\footnote{In general, as will be demonstrated with examples,
our analysis can be carried out for (integral)
input-to-output/state stable systems as well.} \cite{Sontag_1990}
with respect to a compact set $\mathcal{A}$:

\begin{assume}[{\normalfont \it Contracting dynamics}]\label{assume:converging_dynamics}
\begin{equation}\label{eq:attracting}
\mathcal{S}_a:  \ \ \ \norms{\bfx(t)}\leq
\beta(\norms{\bfx(t_0)},t-t_0) + c \|u_a(t)\|_{\infty,[t_0,t]}, \
\forall t_0\in\Real_+, \ t\geq t_0
\end{equation}
where the function $\beta(\cdot,\cdot)\in\mathcal{KL}$, and $c>0$
is some positive constant.
\end{assume}
The function $\beta(\cdot,\cdot)$ in (\ref{eq:attracting})
specifies the contraction property of the unperturbed dynamics of
$\mathcal{S}_a$. In other words it models the rate with which the
system forgets its initial conditions $\bfx_0$, if left
unperturbed. Propagation of the input to output is estimated in
terms of a continuous mapping, $c \|u_a(t)\|_{\infty,[t_0,t]}$,
which, in our case, is chosen for simplicity to be linear. Notice
that this mapping should not necessarily be contracting. In what
follows we will assume that the function $\beta(\cdot,\cdot)$ and
constant $c$ are known or can be estimated a-priori.

For systems $\mathcal{S}_a$, of which a model is given by a system
of ordinary differential equations
\begin{equation}\label{eq:contraction_diff_eq}
\dot{\bfx}=\bff_x(\bfx,u_a), \
\bff_x(\cdot,\cdot)\in\mathcal{C}^1,
\end{equation}
Assumption \ref{assume:converging_dynamics} is equivalent, for
instance, to the combination of the following
properties\footnote{For a comprehensive characterization of the
input-to-state stability and detailed mathematical arguments we
refer to the paper by E.D. Sontag and Y. Wang
\cite{Sontag_1995}.}:
\begin{enumerate}
\item let $u_a(t)\equiv 0$ for all $t$, then set $\mathcal{A}$ is
Lyapunov stable and globally attracting for
(\ref{eq:contraction_diff_eq});

\item for all $u_a\in\mathcal{U}_a$ and $\bfx_0\in\Real^n$ there
exists a non-decreasing function
$\kappa:\Real_+\rightarrow\Real_+: \ \kappa(0)=0$ such that
\[
\inf_{t\in[0,\infty)} \norms{\bfx(t)}\leq
\kappa(\|u_a(t)\|_{\infty,[t_0,\infty)})
\]
\end{enumerate}

The system $\mathcal{S}_w$ stands for the searching or wandering
dynamics. We will consider  $\mathcal{S}_w$ subject to the
following  conditions:

\begin{assume}[{\normalfont \it Wandering dynamics}]\label{assume:searching_dynamics}

The system $\mathcal{S}_w$ is forward-complete:
\[
u_w(t)\in \mathcal{U}_w\Rightarrow \bfz(t)\in\mathcal{Z}, \
\forall \ t\geq t_0, \ t_0\in\Real_+
\]
and there exists an "output" function $h:\Real^m\rightarrow\Real$,
and two "bounding" functions $\gamma_0\in \mathcal{K}_{\infty,e}$,
$\gamma\in\mathcal{K}_{\infty,e}$ such that the following integral
inequality holds:
\begin{equation}\label{eq:integral}
\mathcal{S}_w: \   \int_{t_0}^{t}\gamma_1(u_w(\tau))d\tau \leq
h(\bfz(t_0))-h(\bfz(t))\leq
\int_{t_0}^{t}\gamma_0(u_w(\tau))d\tau, \ \forall \ t\geq t_0, \
t_0\in\Real_+
\end{equation}
\end{assume}
In case system $\mathcal{S}_w$ is specified in terms of
vector-fields
\begin{equation}\label{eq:wandering_diff_eq}
\dot{\bfz}=\bff_z(\bfz,u_w),  \
\bff_z(\cdot,\cdot)\in\mathcal{C}^1,
\end{equation}
Assumption \ref{assume:searching_dynamics} can be viewed, for
example, as postulating the existence of a function
$h:\Real^m\rightarrow\Real_+$ of which the evolution in time is a
mere integration of the input $u_w(t)$. In general, for $u_w: \
u_w(t)\geq 0 \ \ \forall \ t\in\Real_+$, inequality
(\ref{eq:integral}) implies {\it monotonicity} of function
$h(\bfz(t))$ in $t$. Regarding the function $\gamma_0(\cdot)$ in
(\ref{eq:integral}), we assume that for any $M\in\Real_+$ there
exist functions $\gamma_{0,1}:\Real_+\rightarrow\Real_+$ and
${\gamma}_{0,2}:\Real_+\rightarrow\Real_+$, such that
\begin{equation}\label{eq:assume:gamma_0}
\gamma_0(a\cdot b) \leq \gamma_{0,1}(a)\cdot {\gamma}_{0,2}(b),
\forall \ a,b\in[0,M].
\end{equation}
Requirement (\ref{eq:assume:gamma_0}) is a technical assumption
which will be used in the formulation and proof of the main
results of the paper. Yet, it is not too restrictive; it holds,
for instance, for a wide class of locally Lipschitz functions
$\gamma_0(\cdot): \ \gamma_0(a\cdot b)\leq L_{0}(M) \cdot (a\cdot
b)$, $L_{0}(M)\in\Real_+$. Another example for which the
assumption holds is the class of polynomial functions
$\gamma_0(\cdot): \ \gamma_0(a\cdot b)= (a\cdot b)^p= a^p \cdot
b^p$, $p>0$. No further restrictions will be imposed a-priori on
$\mathcal{S}_a$, $\mathcal{S}_w$.

\begin{figure}[t!]
\begin{center}
\hspace{5mm}
\begin{minipage}[h]{0.3\linewidth}
\begin{center}
a
\end{center}
\begin{center}
\includegraphics[width=170pt]{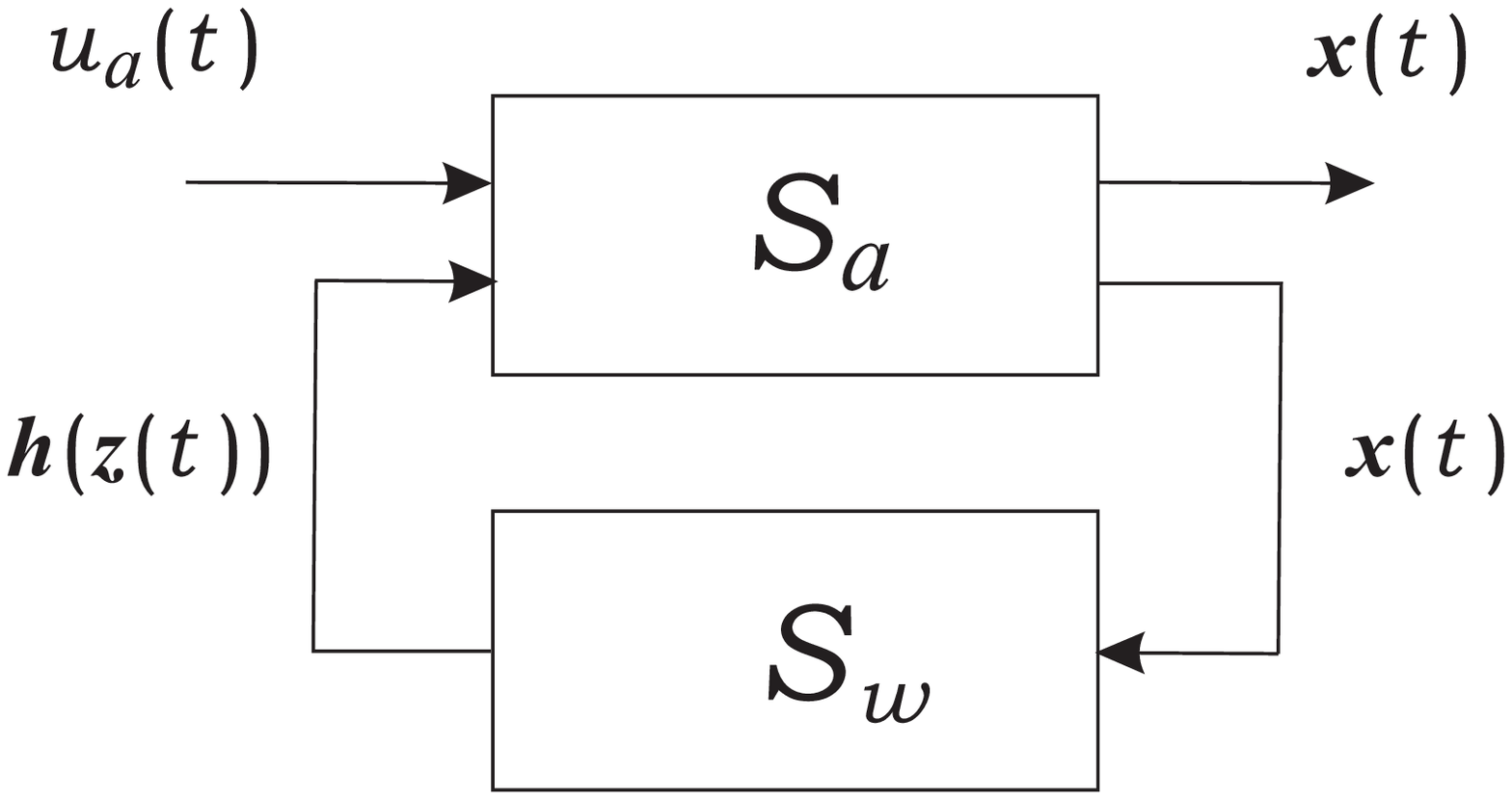}
\end{center}
\end{minipage}
\hspace{17mm}
\begin{minipage}[h]{0.4\linewidth}
\begin{center}
b

\vspace{5mm}

\includegraphics[width=210pt]{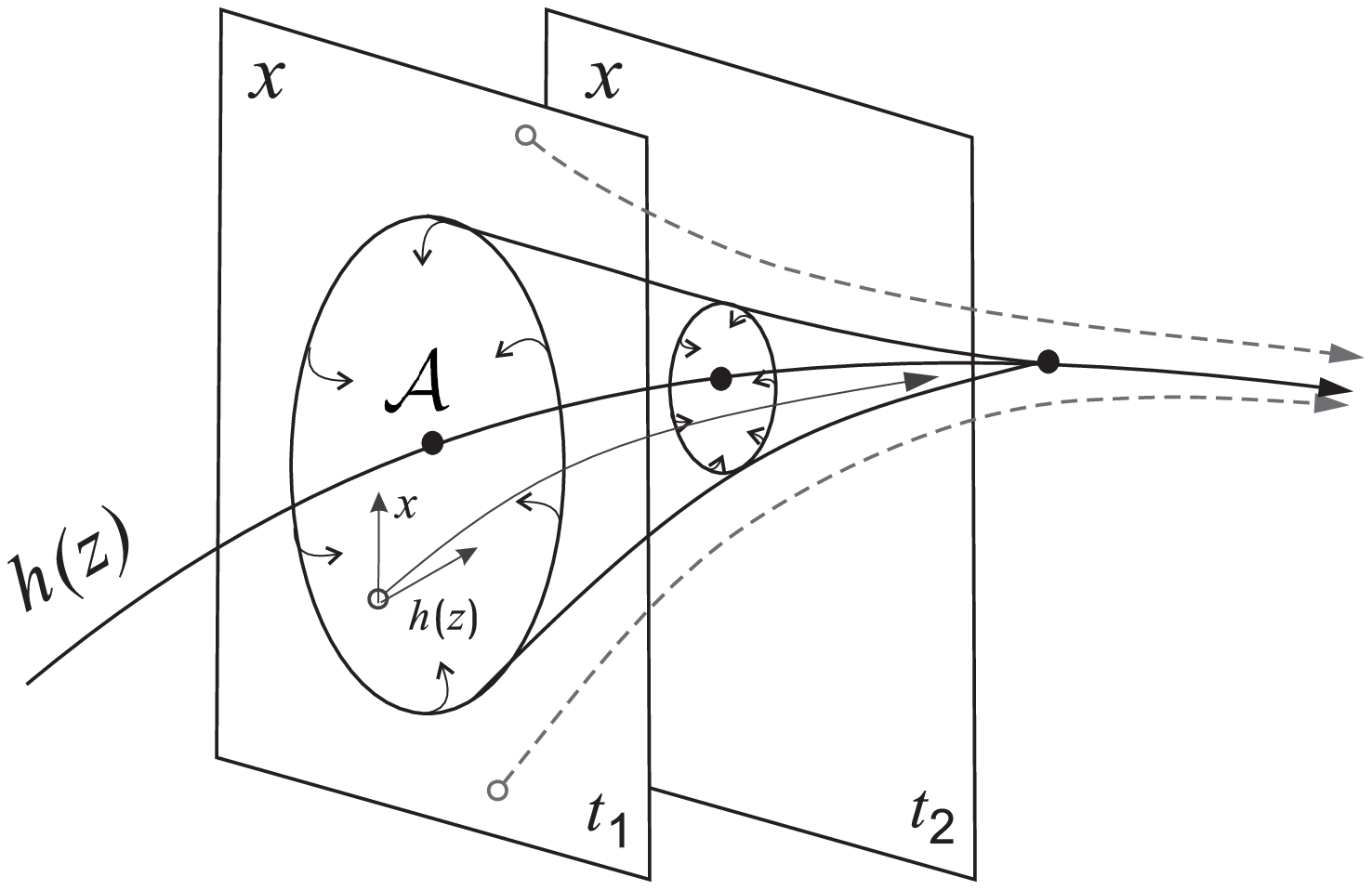}
\end{center}
\end{minipage}
\end{center}
\caption{\small {\bf a.} The class of interconnected systems $S_a$
and $S_w$. System $\mathcal{S}_a$, the ``contracting system", has
an attracting invariant set $\mathcal{A}$ in its state space,
system $\mathcal{S}_w$ does not necessarily have an attracting
set. It represents the ``wandering" dynamics. A typical example of
such behavior is the dynamics of the flow in a neighborhood of a
saddle point in three-dimensional space
(diagram {\bf b}). }\label{fig:SaSw_interconnection}
\end{figure}

Now consider the interconnection of (\ref{eq:attracting}),
(\ref{eq:integral}) with coupling $u_a(t)=h(\bfz(t))$, and
$u_s(t)=\norms{\bfx(t)}$. Equations for the combined system can be
written as:
\begin{equation}\label{eq:interconnection}
\begin{split}
\norms{\bfx(t)}\leq &\beta(\norms{\bfx(t_0)},t-t_0) + c
\|h(\bfz(t))\|_{\infty,[t_0,t]}\\
 \int_{t_0}^{t}\gamma_1(\norms{\bfx(\tau)})d\tau \leq
& h(\bfz(t_0))-h(\bfz(t))\leq
\int_{t_0}^{t}\gamma_0(\norms{\bfx(\tau)})d\tau,
\end{split}
\end{equation}
A diagram illustrating the general structure of the entire system
(\ref{eq:interconnection}) is given in figure
\ref{fig:SaSw_interconnection}.

Equations (\ref{eq:interconnection}) capture the relevant
interplay between contracting, $\mathcal{S}_a$, and wandering,
$\mathcal{S}_w$, dynamics inherent to a variety of searching
strategies in the realm of optimization, (\ref{eq:global_search}),
and interconnections (\ref{eq:cascade_intro}) in general systems
theory. In addition, this kind of interconnection describes the
behavior of systems which undergo transcritical or saddle-node
bifurcations. Consider for instance the following system:
\begin{equation}\label{eq:saddle_node1}
\begin{split}
\dot{x}_1&=-x_1 + x_2\\
\dot{x}_2&=\varepsilon + \gamma x_1^2, \gamma>0
\end{split}
\end{equation}
where the parameter $\varepsilon$ varies from negative to positive
values. At $\varepsilon=0$ stable and unstable equilibria collide
leading to the cascade satisfying equations
(\ref{eq:interconnection}). An alternative bifurcation scenario
could be represented by system:
\begin{equation}\label{eq:saddle_node2}
\begin{split}
\dot{x}_1&=-x_1 + x_2\\
\dot{x}_2&=\varepsilon + \gamma x_2^2, \gamma>0,
\end{split}
\end{equation}
In this case, however, the dynamics of the variable $x_2$ is {\it
independent} of $x_1$, and analysis of asymptotic behavior of
(\ref{eq:saddle_node2}) reduces to the analysis of each equation
separately. Thus systems like (\ref{eq:saddle_node2}) are easier
to deal with than (\ref{eq:saddle_node1}). This constitutes an
additional motivation for the present approach.

When analyzing the asymptotic behavior of interconnection
(\ref{eq:interconnection}) we will address the following set of
questions: is there a set (a weak trapping set in the system state
space) such that the trajectories  which start in this set are
bounded? It is natural to expect that the existence of such a set
depends on the specific  functions $\gamma_0(\cdot)$,
$\gamma_1(\cdot)$ in (\ref{eq:interconnection}), on properties of
$\beta(\cdot,\cdot)$, and on values of $c$. In case such a set
exists and could be defined, the next questions are therefore:
where  will the trajectories converge and how can these these
domains be characterized?




\section{Main Results}

In this section we provide a formal statement of the  main results
of our present study. In Section \ref{subsec:small_gain}, we
formulate conditions ensuring that there exists a point
$\bfx_0\oplus\bfz_0$ such that the $\omega$-limit set of
$\bfx_0\oplus\bfz_0$ is bounded in the following sense:
\begin{equation}\label{eq:omega_boundedness}
\norms{\omega_\bfx(\bfx_0\oplus\bfz_0)}<\infty, \ \
|h(\omega_\bfz(\bfx_0\oplus\bfz_0))|<\infty
\end{equation}
These conditions and also a specification of the set
$\Omega_\gamma$ of points $\bfx'\oplus\bfz'$ for which the
$\omega$-limit set satisfies property (\ref{eq:omega_boundedness})
are provided in Theorem \ref{theorem:non_uniform_small_gain}.

\begin{figure}[t!]
\begin{center}
\includegraphics[width=300pt]{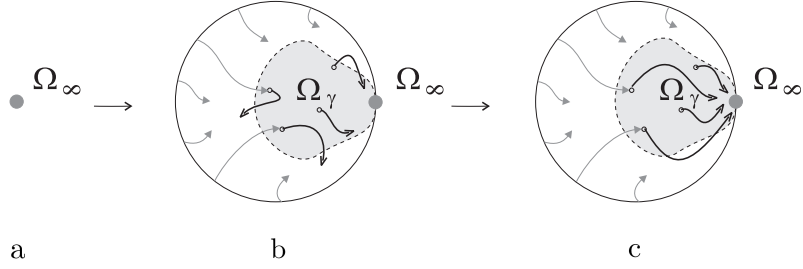}
\end{center}
\caption{\small Emergence of a weak (Milnor) attracting set
$\Omega_\infty$. Panel $a$ depicts the target invariant set
$\Omega_\infty$ as a filled circle. First (Theorem
\ref{theorem:non_uniform_small_gain}), we investigate whether
 a domain $\Omega_\gamma\subset\Real^n\times\Real^m$ exists such that
$\norms{\bfx(t)}$, $h(\bfz(t))$ are bounded for all
$\bfx_0\oplus\bfz_0\in\Omega_\gamma$. In the text we refer to this
set as a weak trapping region or simply a trapping region. The
trapping region is shown as a grey domain in panel $b$. In
principle, the system's states can eventually leave the domain
$\Omega_\gamma$. They must, however, satisfy equation
(\ref{eq:omega_boundedness}), ensuring boundedness of
$\norms{\bfx(t)}$, $h(\bfz(t))$. As a result they will dwell
within the region shown as a circle in panel $b$. Notice that
neither this domain, nor the previous, need be neighborhoods of
$\Omega_\infty$. Second (Lemma's
\ref{lem:steady_state_convergence},
\ref{lem:steady_state_average_convergence}, Corollary
\ref{cor:non_uniform_small_gain_attracting_set}), we provide
conditions which lead to the emergence of a weak attracting set in
the trapping region $\Omega_\gamma$. This is illustrated in panel
$c$.}\label{fig:emergence_of_attractor}
\end{figure}

In order to verify whether an attracting set exists in
$\omega(\Omega_\gamma)$ that is a subset of
$\omega(\Omega_\gamma)$ we use an additional characterization of
the contracting system $\mathcal{S}_a$. In particular,  we
introduce the intuitively clear notion of the input-to state {\it
steady-state characteristics}\footnote{A more precise definition
of the steady-state characteristics is given in Section
\ref{subsec:attracting_sets}} of a system. It is possible to show
that in case system $\mathcal{S}_a$ has a steady-state
characteristic, then there exists an attracting set in
$\omega(\Omega_\gamma)$ and this set is uniquely defined by the
zeros of the steady-state characteristics of $\mathcal{S}_a$.
A diagram illustrating the steps of our analysis is provided in
Fig. \ref{fig:emergence_of_attractor}, as well as the sequence of
conditions leading to the emergence of the attracting set in
(\ref{eq:interconnection}).


\subsection{Emergence of the trapping region. Small-gain
conditions}\label{subsec:small_gain}

Before we formulate the main results of this subsection let us
first comment briefly on the machinery of our analysis. First of
all we introduce three sequences
\[
\mathcal{S}=\{\sigma_i\}_{i=0}^\infty,  \ \ \
\Xi=\{\xi_i\}_{i=0}^\infty, \ \
\mathcal{T}=\{\tau_i\}_{i=0}^\infty
\]
The first sequence, $\mathcal{S}$, partitions the interval
$[0,h(\bfz_0)]$, $h(\bfz_0)>0$ into the union of shrinking
subintervals $H_i$:
\begin{equation}\label{eq:partitioned_interval}
[0,h(\bfz_0)]=\cup_{i=0}^\infty H_i, \ H_i=
[\sigma_{i}h(\bfz_0),\sigma_{i+1}h(\bfz_0)]
\end{equation}
For the sake of transparency, let us define this property formally
in the form of Condition \ref{assume:partition_of_z}
\begin{cond}[Partition of $\bfz_0$]\label{assume:partition_of_z} The sequence $\mathcal{S}$
is strictly monotone and converging
\begin{equation}\label{eq:partitioning_sequence}
\{\sigma_n\}_{n=0}^\infty: \ \lim_{n\rightarrow\infty}\sigma_n=0,
\ \sigma_0=1
\end{equation}
\end{cond}
Sequences $\Xi$ and $\mathcal{T}$ will specify the desired rates
$\xi_i\in\Xi$ of the contracting dynamics (\ref{eq:attracting}) in
terms of function $\beta(\cdot,\cdot)$  and  time
$T_i>\tau_i\in\mathcal{T}$. Let us, therefore, impose the
following constraint on the choice of $\Xi$, $\mathcal{T}$.
\begin{cond}[Rate of contraction, Part 1]\label{assume:desired_contraction_rates_1} For the given sequences $\Xi$,
$\mathcal{T}$ and function $\beta(\cdot,\cdot)\in\mathcal{KL}$ in
(\ref{eq:attracting}) the following inequality holds:
\begin{equation}\label{eq:tau_system}
\beta(\cdot,T_i)\leq \xi_i\beta(\cdot,0), \ \forall  \
T_i\geq\tau_i
\end{equation}
\end{cond}
Condition \ref{assume:desired_contraction_rates_1} states that for
the given, yet arbitrary, factor $\xi_i$ and time instant $t_0$,
the amount of time $\tau_i$ is needed  for the state $\bfx$ in
order to reach the domain:
\[
\norms{\bfx}\leq \xi_i \beta(\norms{\bfx(t_0)},0)
\]
In order to specify the desired convergence rates $\xi_i$,  it
will be necessary to define another measure in addition to
(\ref{eq:tau_system}). This is a measure of the propagation of
initial conditions $\bfx_0$ and input $h(\bfz_0)$ to the state
$\bfx(t)$ of the contracting dynamics (\ref{eq:attracting}) when
the system travels in $h(\bfz(t))\in [0,h(\bfz_0)]$. For this
reason we introduce two systems of functions, $\Phi$ and
$\Upsilon$:
\begin{equation}\label{eq:phi_system}
\Phi: \
\begin{array}{ll}
\phi_j(s)&=\phi_{j-1}\circ\rho_{\phi,j}(\xi_{i-j}\cdot\beta(s,0)),
\ j=1,\dots,i\\
\phi_0(s)&=\beta(s,0)
\end{array}
\end{equation}
\begin{equation}\label{eq:u_system}
\Upsilon: \
\begin{array}{ll}
\upsilon_j(s)&=\phi_{j-1}\circ\rho_{\upsilon,j}(s),
\ j=1,\dots,i\\
\upsilon_0(s)&=\beta(s,0)
\end{array}
\end{equation}
where the functions $\rho_{\phi,j}, \
\rho_{\upsilon,j}\in\mathcal{K}$ satisfy the following inequality
\begin{equation}\label{eq:rho_system}
\phi_{j-1}(a+b)\leq
\phi_{j-1}\circ\rho_{\phi,j}(a)+\phi_{j-1}\circ\rho_{\upsilon,j}(b)
\end{equation}
Notice that in case $\beta(\cdot,0)\in\mathcal{K}_\infty$ the
functions $\rho_{\phi,j}(\cdot)$, $\rho_{\upsilon,j}(\cdot)$ will
always exist \cite{Jiang_1994}. The properties of sequence $\Xi$
which ensure the desired propagation rate of the influence of
initial condition $\bfx_0$ and input $h(\bfz_0)$ to the state
$\bfx(t)$ are specified in Condition \ref{assume:phi_system}.
\begin{cond}[Rate of contraction, Part 2]\label{assume:phi_system} The sequences
\[
\sigma_n^{-1}\cdot \phi_n(\norms{\bfx_0}), \ \
\sigma_n^{-1}\cdot\left(\sum_{i=0}^n\upsilon_i(c|h(\bfz_0)|\sigma_{n-i})\right),
\ n=0,\dots,\infty
\]
are bounded from above, e.g. there exist functions
$B_1(\|\bfx_0\|)$, $B_2(|h(\bfz_0)|,c)$ such that
\begin{equation}\label{eq:contraction_1}
\sigma_n^{-1}\cdot \phi_n(\norms{\bfx_0})\leq B_1(\norms{\bfx_0})
\end{equation}
\begin{equation}\label{eq:contraction_2}
\sigma_n^{-1}\cdot\left(\sum_{i=0}^n\upsilon_i(c|h(\bfz_0)|\sigma_{n-i})\right)\leq
B_2(|h(\bfz_0)|,c)
\end{equation}
for all $n=0,1,\dots,\infty$
\end{cond}
For a large class of functions $\beta(s,0)$, for instance those
that are Lipschitz in $s$, these conditions reduce to more
transparent ones which can always be satisfied by an appropriate
choice of sequences $\Xi$ and $\mathcal{S}$. This case is
considered in detail as a corollary of our main results in section
\ref{subsec:separable_dynamics}.

\begin{figure}[t!]
\begin{tabular*}{\linewidth}[t]{|c|c|}\hline
 & \\
Standard & Proposed\\
& \\
\hline & \\
\begin{minipage}[c]{0.45\linewidth}
\begin{enumerate}
\item[1) ] Domain of attraction is a neighborhood
\end{enumerate}
\end{minipage}
&
\begin{minipage}[c]{0.45\linewidth}
\begin{enumerate}
 \item[1) ] Domain of attraction is a set of positive measure (not necessarily a neighborhood)
\end{enumerate}
\end{minipage}\\
& \\
\begin{minipage}[c]{0.45\linewidth}
\begin{enumerate}
\item[2) ] Implies Lyapunov stability
\end{enumerate}
\end{minipage}
&
\begin{minipage}[c]{0.45\linewidth}
\begin{enumerate}
 \item[2) ] Allows to analyze convergence in Lyapunov-unstable
 systems
\end{enumerate}
\end{minipage}\\
& \\
\begin{minipage}[c]{0.475\linewidth}
\begin{center}
\includegraphics[width=160pt]{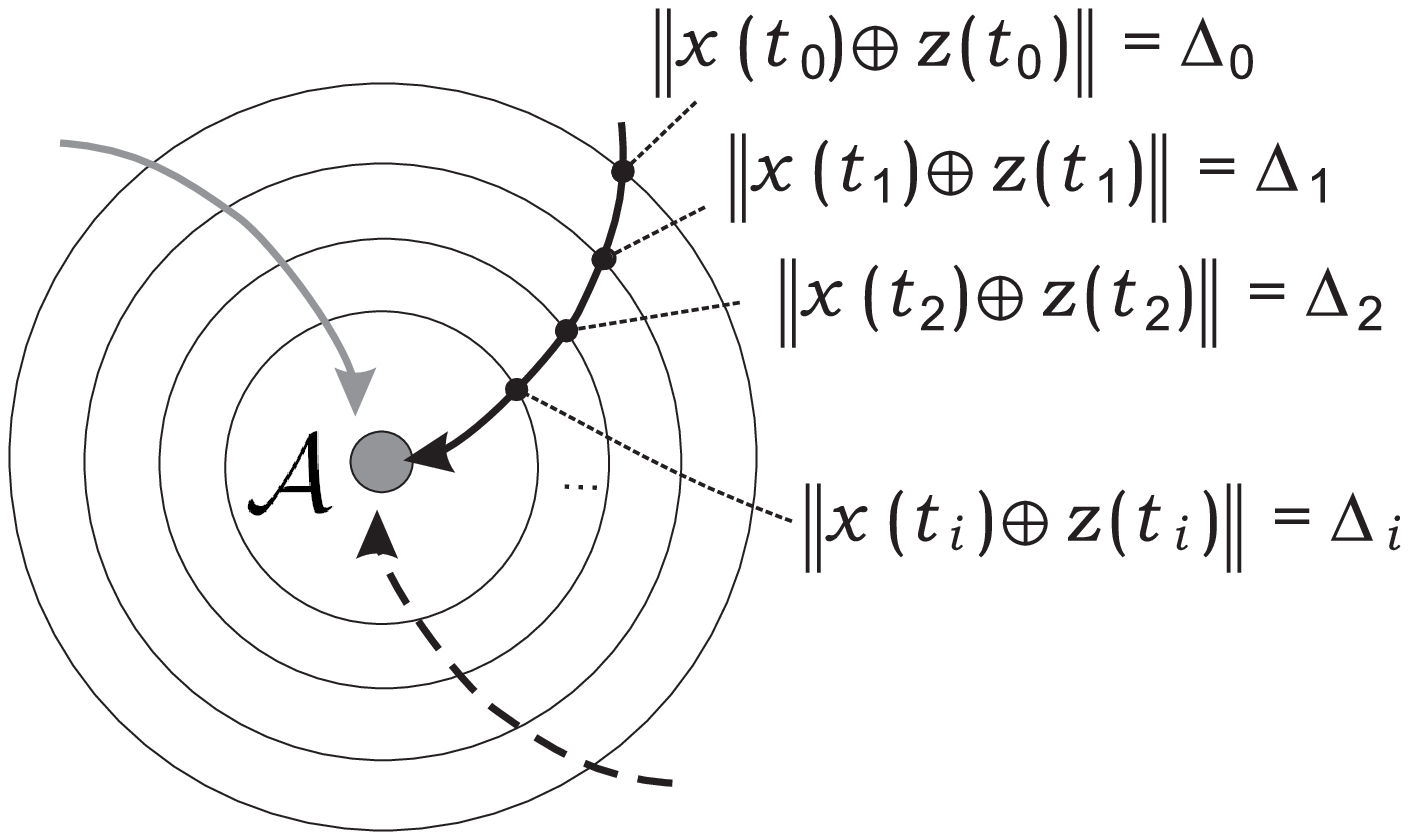}
\end{center}
\end{minipage}
&
\begin{minipage}[c]{0.475\linewidth}
\begin{center}
\includegraphics[width=180pt]{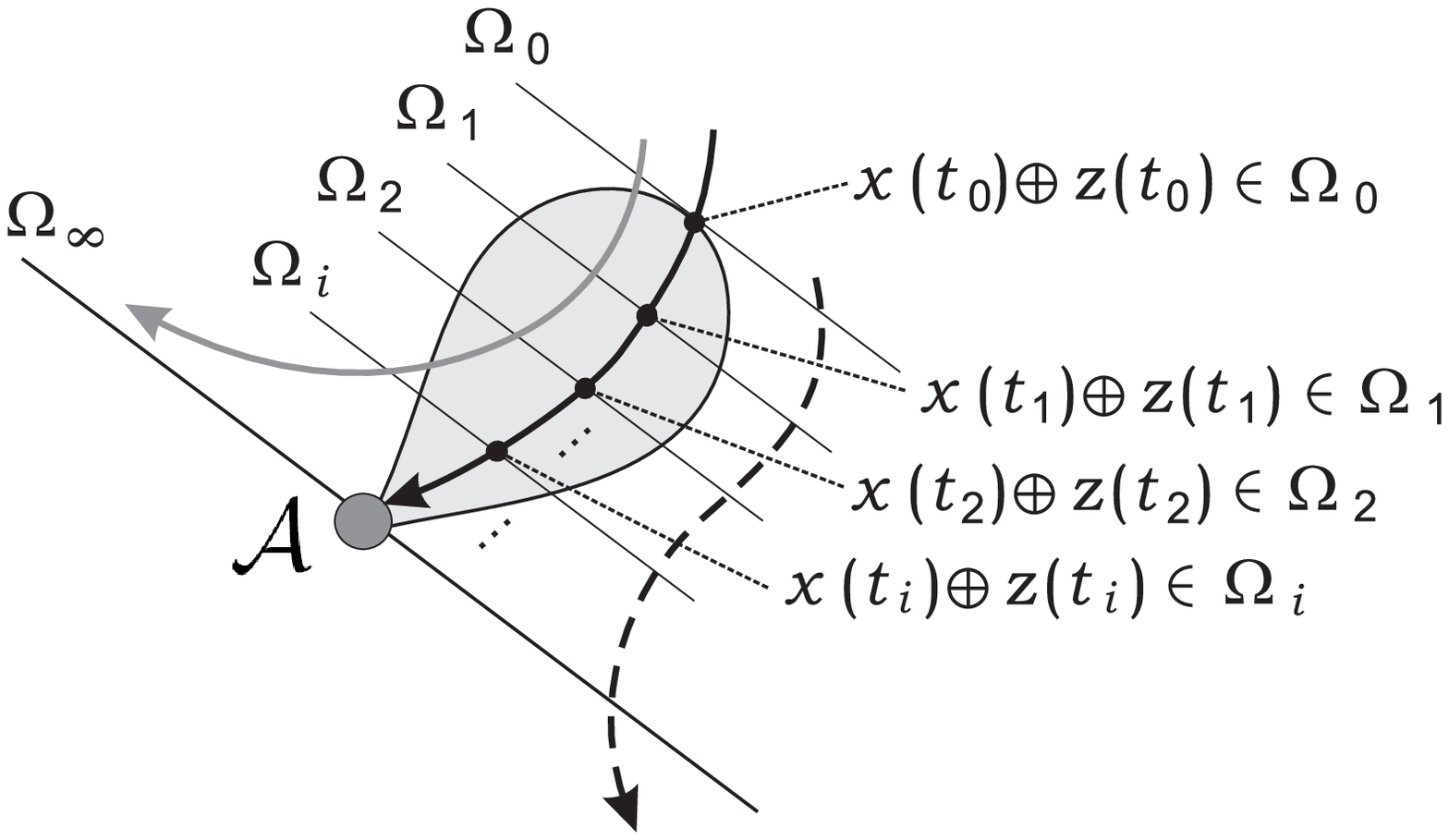}
\end{center}
\end{minipage}
\\
& \\
\hline & \\
\begin{minipage}[c]{0.45\linewidth}
Given: a sequence of diverging time instances $t_i$
\end{minipage}
&
\begin{minipage}[c]{0.45\linewidth}
Given: a sequence of sets $\Omega_i$ whose distance $\Delta_i$ to
$\mathcal{A}$ is converging to zero
\end{minipage}\\
& \\
\hline & \\
\begin{minipage}[c]{0.45\linewidth}
Prove: convergence of norms
$\|\bfx(t_i)\oplus\bfz(t_i)\|=\Delta_i$ to zero
\end{minipage}
&
\begin{minipage}[c]{0.45\linewidth}
Prove: divergence of $\{t_i\}$, where $t_i: \
\bfx(t_i)\oplus\bfz(t_i)\in\Omega_i$
\end{minipage}\\
& \\
\hline
\end{tabular*}
\caption{\small Key differences between the conventional concept
of convergence (left panel)  and the concept of weak, non-uniform,
convergence (right panel). In the uniform case, trajectories which
start in a neighborhood of $\mathcal{A}$ remain in a neighborhood
of $\mathcal{A}$ (solid and dashed lines). In the non-uniform
case, only a fraction of the initial conditions in a neighborhood
of $\mathcal{A}$ will produce trajectories which remain in a
neighborhood of $\mathcal{A}$ (solid black line). In the most
general case a necessary condition for this to happen is that the
sequence $\{t_i\}$ diverges. In our current problem statement
divergence of $\{t_i\}$ implies boundedness of $\norms{\bfx(t)}$.
To show state boundedness and convergence of $\bfx(t)$ to
$\mathcal{A}$ an additional information on the system dynamics
will be required.}\label{table:non_uniform_attractivity}
\end{figure}

The main differences  between the standard and the presently
proposed approaches  for the analysis of asymptotic behavior of
dynamical systems are illustrated with figure
\ref{table:non_uniform_attractivity}.
In order to prove the emergence of the trapping region we consider
the following collection of volumes induced by the sequence
$\mathcal{S}_i$  and the corresponding partition
(\ref{eq:partitioned_interval}) of the interval $[0, h(\bfz_0)]$:
\begin{equation}\label{eq:shrinking_volumes}
\Omega_i=\{\bfx\in\mathcal{X},\bfz\in\mathcal{Z}| \ h(\bfz(t))\in
H_i\}
\end{equation}
For the given initial conditions $\bfx_0\in\mathcal{X}$,
$\bfz_0\in\mathcal{Z}$ two alternative possibilities exist. First,
the trajectory $\bfx(t,\bfx_0)\oplus\bfz(t,\bfz_0)$  stays in some
$\Omega'\subset\Omega_0$ for all $t>t'$, $t'\geq t_0$. Hence for
$t\rightarrow\infty$ the state will converge into
\begin{equation}\label{eq:domain_of_convergence}
\Omega_a=\{\bfx\in\mathcal{X}, \ \bfz\in\mathcal{Z}|
\norms{\bfx}\leq c\cdot h(\bfz_0), \ \bfz: \ h(\bfz)\in[0,
h(\bfz_0)]\}
\end{equation}
Second, the trajectory $\bfx(t,\bfx_0)\oplus\bfz(t,\bfz_0)$
subsequently enters the volumes $\Omega_j$, and $t_j$ are the time
instances when it hits the hyper-surfaces
$h(\bfz(t))=h(\bfz_0)\sigma_j$. Then the state of the coupled
system stays in $\Omega_0$ only if the sequence
$\{t_i\}_{i=0}^\infty$ diverges. Theorem
\ref{theorem:non_uniform_small_gain} provides the conditions
specifying the latter case in terms of properties of sequences
$\mathcal{S}$, $\Xi$, $\mathcal{T}$ and function $\gamma_0(\cdot)$
in (\ref{eq:interconnection}).

\begin{thm}[Non-uniform Small-gain Theorem]\label{theorem:non_uniform_small_gain} Let systems $\mathcal{S}_a$, $\mathcal{S}_w$  be given and satisfy Assumptions \ref{assume:converging_dynamics},
\ref{assume:searching_dynamics}. Consider their interconnection
(\ref{eq:interconnection}) and suppose there exist sequences
$\mathcal{S}$, $\Xi$, and $\mathcal{T}$ satisfying Conditions
\ref{assume:partition_of_z}--\ref{assume:phi_system}. In addition,
suppose that the following conditions hold:

1) There exists a positive number $\Delta_0>0$ such that
\begin{equation}\label{eq:tau_system_1}
\frac{1}{\tau_i}\frac{(\sigma_i-\sigma_{i+1})}{\gamma_{0,1}(\sigma_i)}\geq
\Delta_0 \ \forall  \ i=0,1,\dots,\infty
\end{equation}

2) The set $\Omega_\gamma$ of all points $\bfx_0$, $\bfz_0$
satisfying the inequality
\begin{equation}\label{eq:non_uniform_small_gain}
\gamma_{0,2}(B_1(\norms{\bfx_0})+B_2(|h(\bfz_0)|,c)+c
|h(\bfz_0)|)\leq h(\bfz_0)\Delta_0
\end{equation}
is not empty.

3) Partial sums of  elements from  $\mathcal{T}$ diverge:
\begin{equation}\label{eq:tau_system_2}
\sum_{i=0}^\infty\tau_i=\infty
\end{equation}

\noindent Then for all $\bfx_0$, $\bfz_0\in\Omega_\gamma$ the
state $\bfx(t,\bfz_0)\oplus\bfz(t,\bfz_0)$ of system
(\ref{eq:interconnection}) converges into the set specified by
(\ref{eq:domain_of_convergence})
\begin{equation}
\begin{split}
\Omega_a=\{\bfx\in\mathcal{X}, \ \bfz\in\mathcal{Z}|
\norms{\bfx}\leq c\cdot h(\bfz_0), \ \bfz: \ h(\bfz)\in[0,
h(\bfz_0)]\}\nonumber
\end{split}
\end{equation}
\end{thm}
The proof of the theorem is provided in Appendix 1.

The major difference between the conditions of Theorem
\ref{theorem:non_uniform_small_gain} and those of conventional
small-gain theorems \cite{Zames:66},\cite{Jiang_1994} is that the
latter involve only input-output or input-state mappings.
Formulating conditions for state boundedness of the
interconnection in terms of input-output or input-state mappings
is possible in the traditional case because the interconnected
systems are assumed to be input-to-state stable. Hence their
internal dynamics can be neglected. In our case, however, the
dynamics of $\mathcal{S}_w$ is generally unstable in the Lyapunov
sense. Hence, in order to ensure boundedness of $\bfx(t,\bfx_0)$
and $h(\bfz(t,\bfz_0))$, the rate/degree of stability of
$\mathcal{S}_a$ should be taken into account. Roughly speaking,
system $\mathcal{S}_a$ should ensure a sufficiently high degree of
contraction in $\bfx_0$ while the input-output response of
$\mathcal{S}_w$ should be sufficiently small. The rate of
contraction in $\bfx_0$ of $\mathcal{S}_a$, according to
(\ref{eq:attracting}), is specified in terms of the function
$\beta(\cdot,\cdot)$. Properties of this function that are
relevant for convergence  are explicitly accounted for in
Condition \ref{assume:phi_system} and (\ref{eq:tau_system_2}). The
domain of admissible initial conditions and actually the
small-gain condition (input-state-output properties of
$\mathcal{S}_w$ and $\mathcal{S}_a$) are defined by
(\ref{eq:tau_system_1}), (\ref{eq:non_uniform_small_gain})
respectively. Notice also that $\Omega_\gamma$ is not necessarily
a neighborhood of ${\Omega}_a$, thus the convergence ensured by
Theorem \ref{theorem:non_uniform_small_gain} is allowed to
be non-uniform in $\bfx_0$, $\bfz_0$.

\subsection{Characterization of the attracting
set}\label{subsec:attracting_sets}

Even for interconnections of  Lyapunov-stable systems, small-gain
conditions usually are effective  merely for establishing
boundedness of states or outputs. Yet, even in the setting of
Theorem \ref{theorem:non_uniform_small_gain} it is still possible
to derive estimates (such as, for instance
(\ref{eq:domain_of_convergence})) of the domains to which the
state will converge. These estimates, however, are often too
conservative. If a more precise characterization of these domains
is required, additional information on the dynamics of systems
$\mathcal{S}_a$ and $\mathcal{S}_w$ will be needed. The question,
therefore, is how detailed this information should be? It appears
that some additional knowledge of the steady-state characteristics
of system $\mathcal{S}_a$ is sufficient to improve the estimates
(\ref{eq:domain_of_convergence}) substantially.

Let us formally introduce the notion of steady-state
characteristic as follows:
%

\begin{defn}\label{defn:steady_state_norm} We say that system
(\ref{eq:attracting}) has steady-state characteristic
$\chi:\Real\rightarrow\mathcal{S}\{\Real_+\}$ with respect to the
norm $\norms{\bfx}$ if and only if for each constant $\bar{u}_a$
the following holds:
\begin{equation}\label{eq:steady_state}
\forall \ u_a(t)\in \mathcal{U}_a: \
\lim_{t\rightarrow\infty}u_a(t)=\bar{u}_a \ \Rightarrow \
\lim_{t\rightarrow\infty}\norms{\bfx(t)}\in\chi(\bar{u}_a)
\end{equation}
\end{defn}
The key property captured by Definition
\ref{defn:steady_state_norm} is that there exists a limit of
$\norms{\bfx(t)}$ as $t\rightarrow\infty$, provided that the limit
for $u_a(t)$, $t\rightarrow\infty$ is defined and constant. Notice
that the mapping $\chi$ is set-valued. This means that for each
$\bar{u}_a$ there is a set $\chi(\bar{u}_a)\subset\Real_+$ such
that $\norms{\bfx(t)}$ converges to an element of
$\chi(\bar{u}_a)$ as $t\rightarrow\infty$. Therefore, our
definition allows a fairly large amount of uncertainty for
$\mathcal{S}_a$. It will be of essential importance, however, that
such characterization exists for the system $\mathcal{S}_a$.

Clearly, not every system obeys a steady-state characteristic
$\chi(\cdot)$ of Definition \ref{defn:steady_state_norm}. There
are relatively simple systems of which the state does not converge
even in the "norm" sense for constant converging inputs (condition
(\ref{eq:steady_state})). In mechanics, physics, and biology such
systems encompass the large class of nonlinear oscillators which
can be excited by constant inputs. In order to take such systems
into consideration, we introduce a weaker notion, that of
steady-state characteristic {\it on average}:

\begin{defn}\label{defn:steady_state average_norm} We say that system
(\ref{eq:attracting}) has steady-state characteristic on average
$\chi_T:\Real\rightarrow\mathcal{S}\{\Real_+\}$ with respect to
the norm $\norms{\bfx}$ if and only if for each constant
$\bar{u}_a$ and some $T>0$ the following holds:
\begin{equation}\label{eq:steady_state_average}
\forall \ u_a(t)\in \mathcal{U}_a: \
\lim_{t\rightarrow\infty}u_a(t)=\bar{u}_a \ \Rightarrow \
\lim_{t\rightarrow\infty}\int_{t}^{t+T}\norms{\bfx(\tau)}d\tau\in\chi_T(\bar{u}_a)
\end{equation}
\end{defn}

\noindent Steady-state characterizations of  system
$\mathcal{S}_a$ allow  to further specify the asymptotic behavior
of interconnection (\ref{eq:interconnection}). These results are
summarized in Lemmas \ref{lem:steady_state_convergence} and
\ref{lem:steady_state_average_convergence} below.

\begin{lem}\label{lem:steady_state_convergence} Let system
(\ref{eq:interconnection}) be given and $h(\bfz(t,\bfz_0))$ be
bounded for some $\bfx_0,\bfz_0$. Let, furthermore, system
(\ref{eq:attracting}) have steady-state characteristic
$\chi(\cdot): \Real\rightarrow\mathcal{S}\{\Real_+\}$. Then the
following limiting relations hold\footnote{The symbol
$\chi^{-1}(0)$ in equation (\ref{eq:converge_steady_state})
denotes the set: $\chi^{-1}(0)=\bigcup_{\bar{u}_a\in\Real_+}
\bar{u}_a: \ \chi(\bar{u}_a)\ni 0$.}
\begin{equation}\label{eq:converge_steady_state}
\lim_{t\rightarrow\infty} \norms{\bfx(t,\bfx_0)}=0, \
\lim_{t\rightarrow\infty} h(\bfz(t,\bfz_0))\in\chi^{-1}(0)
\end{equation}
\end{lem}

\noindent As follows from Lemma
\ref{lem:steady_state_convergence}, in case the steady-state
characteristic of $\mathcal{S}_a$ is defined, the asymptotic
behavior of interconnection (\ref{eq:interconnection}) is
characterized by the zeroes of the steady-state mapping
$\chi(\cdot)$. For the steady-state characteristics on average a
slightly modified conclusion can be derived.

\begin{lem}\label{lem:steady_state_average_convergence} Let system
(\ref{eq:interconnection}) be given, $h(\bfz(t,\bfz_0))$ be
bounded for some $\bfx_0,\bfz_0$,
$h(\bfz(t,\bfz_0))\in[0,h(\bfz_0)]$, and system
(\ref{eq:attracting}) have steady-state characteristic
$\chi_T(\cdot): \Real\rightarrow\mathcal{S}\{\Real_+\}$ on
average. Furthermore, let there exist a positive constant
$\bar{\gamma}$ such that the
 function $\gamma_1(\cdot)$ in (\ref{eq:integral}) satisfies the
 following constraint:
\begin{equation}\label{eq:converge_steady_state_average_2}
\gamma_1(s)\geq \bar{\gamma} \cdot s, \ \forall s\in[0,\bar{s}],
\bar{s}\in\Real_+: \ \bar{s}>c\cdot h(\bfz_0),
\end{equation}
In addition, suppose that $\chi_T(\cdot)$ has no zeros in the
positive domain, i.e. $0\notin \chi_T(\bar{u}_a)$ for all
$\bar{u}_a>0$. Then
\begin{equation}\label{eq:converge_steady_state_average}
\lim_{t\rightarrow\infty} \norms{\bfx(t,\bfx_0)}=0, \
\lim_{t\rightarrow\infty} h(\bfz(t,\bfz_0))=0
\end{equation}
\end{lem}

An immediate outcome of Lemmas \ref{lem:steady_state_convergence}
and \ref{lem:steady_state_average_convergence} is that in case the
conditions of Theorem \ref{theorem:non_uniform_small_gain} are
satisfied and system (\ref{eq:attracting}) has steady-state
characteristics $\chi(\cdot)$ or $\chi_T(\cdot)$ the domain of
convergence $\Omega_a$ becomes
\begin{equation}\label{eq:new_domain_of_convergence}
\Omega_a=\{\bfx\in\mathcal{X}, \ \bfz\in\mathcal{Z}|
\norms{\bfx}=0, \ \bfz: \ h(\bfz)\in[0, h(\bfz_0)]\}
\end{equation}
It is possible, however, to improve estimate
(\ref{eq:new_domain_of_convergence}) further under additional
hypotheses on system $\mathcal{S}_a$ and $\mathcal{S}_w$ dynamics.
This result is formulated in the corollary below.

\begin{cor}\label{cor:non_uniform_small_gain_attracting_set} Let
system (\ref{eq:interconnection}) be given and satisfy the
assumptions of Theorem \ref{theorem:non_uniform_small_gain}. Let,
in addition,

C1) the flow $\bfx(t,\bfx_0)\oplus\bfz(t,\bfz_0)$ be generated by
a system of autonomous differential equations with locally
Lipschitz right-hand side;

C2) subsystem $\mathcal{S}_w$ be practically
integral-input-to-state stable:
\begin{equation}\label{eq:integral_practically_ISS}
\|\bfz(\tau)\|_{\infty,[t_0,t]}\leq C_z +
\int_{0}^t\gamma_{1}(u_w(\tau))d\tau
\end{equation}
and let function $h(\cdot)\in\mathcal{C}^0$ in (\ref{eq:integral})

C3) system $\mathcal{S}_a$ have steady-state characteristic
$\chi(\cdot)$.

\noindent Then for all $\bfx_0, \bfz_0\in\Omega_\gamma$ the state
of the interconnection converges to the set
\begin{equation}\label{eq:converge_steady_state_set}
\Omega_a=\{\bfx\in\mathcal{X},\bfz\in\mathcal{Z} | \
\norms{\bfx}=0, \ h(\bfz)\in\chi^{-1}(0)\}
\end{equation}
\end{cor}

\begin{figure}
\begin{center}
\includegraphics[width=450pt]{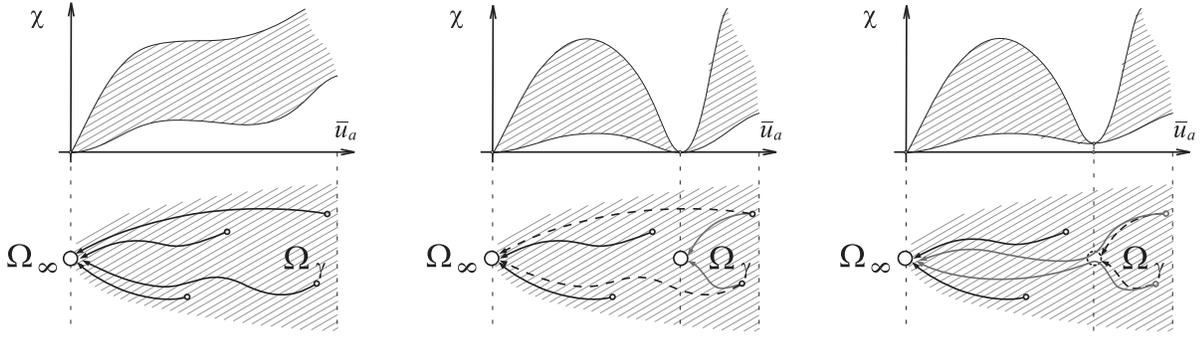}
\end{center}
\caption{Control of the attracting set by means of the system's
steady-state characteristics}\label{fig:attracting_sets}
\end{figure}

As follows from Corollary
\ref{cor:non_uniform_small_gain_attracting_set} zeros of the
steady state characteristic of system $\mathcal{S}_a$ actually
"controls" the domains to which the state of interconnection
(\ref{eq:interconnection}) might potentially converge. This is
illustrated in Fig. \ref{fig:attracting_sets}. Notice also that in
case condition C3 in Corollary
\ref{cor:non_uniform_small_gain_attracting_set} is replaced with
the alternative:

{\it C3') system $\mathcal{S}_a$ has a steady-state characteristic
on average $\chi_T(\cdot)$},

\noindent then it is possible to show that the state converges to
\begin{equation}\label{eq:converge_steady_state_set_average}
\Omega_a=\{\bfx\in\mathcal{X},\bfz\in\mathcal{Z} | \
\norms{\bfx}=0, \ h(\bfz)=0\}
\end{equation}
The proof follows straightforwardly from the proof of Corollary
\ref{cor:non_uniform_small_gain_attracting_set} and  is therefore
omitted.

\subsection{Systems with contracting
dynamics separable in space-time
}\label{subsec:separable_dynamics}

In the previous sections we have presented convergence tests and
estimates of the trapping region, and also characterized  the
attracting sets of interconnection (\ref{eq:interconnection})
under assumptions of uniform asymptotic stability of
$\mathcal{S}_a$ and input-output properties (\ref{eq:integral}),
(\ref{eq:integral_practically_ISS}) of system $\mathcal{S}_w$. The
conditions are given for rather general functions
$\beta(\cdot,\cdot)\in\mathcal{KL}$ in (\ref{eq:attracting}) and
$\gamma_0(\cdot)$, $\gamma_1(\cdot)$ in (\ref{eq:integral}). It
appears, however, that these conditions can be substantially
simplified if  additional properties of $\beta(\cdot,\cdot)$ and
$\gamma_0(\cdot)$ are available. This information is, in
particular, the separability of function $\beta(\cdot,\cdot)$ or,
equivalently, the possibility of factorization:
\begin{equation}\label{eq:separable:dynamics}
\beta(\norms{\bfx},t) \leq \beta_x({\norms{\bfx}})\cdot
\beta_t(t),
\end{equation}
where $\beta_x(\cdot)\in\mathcal{K}$ and
$\beta_t(\cdot)\in\mathcal{C}^0$ is strictly
decreasing\footnote{If $\beta_t(\cdot)$ is not strictly monotone,
it can always be majorized by a strictly decreasing function} with
\begin{equation}\label{eq:separable:dynamics:time}
\lim_{t\rightarrow\infty}\beta_t(t)=0
\end{equation}
In principle, as shown in \cite{Grune:1999}, factorization
(\ref{eq:separable:dynamics}) is achievable for a large class of
uniformly asymptotically stable systems under an appropriate
coordinate transformation. An immediate consequence of
factorization (\ref{eq:separable:dynamics}) is that the elements
of sequence $\Xi$ in Condition
\ref{assume:desired_contraction_rates_1}  are independent of
$\norms{\bfx(t_i)}$. As a result, verification of Conditions
\ref{assume:desired_contraction_rates_1}, \ref{assume:phi_system}
becomes easier. The most interesting case, however, occurs when
the function $\beta_x(\cdot)$ in  the factorization
(\ref{eq:separable:dynamics}) is Lipschitz.  For this class of
functions the  conditions of Theorem
\ref{theorem:non_uniform_small_gain} reduce to a single and easily
verifiable inequality. Let us consider this case in detail.

Without loss of generality, we assume that the state $\bfx(t)$ of
system $\mathcal{S}_a$ satisfies the following equation
\begin{equation}\label{eq:uniformly_attractive_dynamics}
\norms{\bfx(t)}\leq \norms{\bfx(t_0)}\cdot
\beta_t(t-t_0)+c\cdot\|h(\bfz(\tau,\bfz_0))\|_{\infty,[t_0,t]},
\end{equation}
where $\beta_t(0)$ is greater or equal to one. Given that
$\beta_t(t)$ is strictly decreasing, the mapping
$\beta_t:[0,\infty]\mapsto[0,\beta_t(0)]$ is injective. Moreover
$\beta_t(t)$ is continuous, then it is surjective and, therefore,
bijective. In the other words there is a (continuous) mapping
$\beta_t^{-1}:[0,\beta_t(0)]\mapsto\Real_+$:
\begin{equation}\label{eq:beta_inverse}
\beta_{t}^{-1}\circ\beta_t(t)=t, \ \forall \ t>0
\end{equation}
Conditions for emergence of the trapping region for
interconnection (\ref{eq:interconnection}) with dynamics of system
$\mathcal{S}_a$ governed by equation
(\ref{eq:uniformly_attractive_dynamics}) are summarized below:

\begin{cor}\label{cor:non_uniform_small_gain_GAS} Let the interconnection (\ref{eq:interconnection}) be
given,  system $\mathcal{S}_a$ satisfy
(\ref{eq:uniformly_attractive_dynamics}) and function
$\gamma_0(\cdot)$ in  (\ref{eq:integral}) be Lipschitz:
\begin{equation}\label{eq:Lipschitz_gamma}
|\gamma_0(s)|\leq D_{\gamma,0} \cdot | s |
\end{equation}
and domain
\begin{equation}\label{eq:domain_of_convergence_separable}
\Omega_\gamma: \ D_{\gamma,0}\leq
\left(\beta_t^{-1}\left(\frac{d}{\kappa}\right)\right)^{-1}\frac{\kappa-1}{\kappa}
\frac{h(\bfz_0)}{\beta_t(0)\norms{\bfx_0}+\beta_t(0)\cdot c\cdot
|h(\bfz_0)|\left(1+\frac{\kappa}{1-d}\right)+c|h(\bfz_0)|}
\end{equation}
is not empty for some $d<1$, $\kappa>1$. Then for all initial
conditions $\bfx_0$ $\bfz_0\in\Omega_\gamma$  the state
$\bfx(t,\bfx_0)\oplus\bfz(t,\bfz_0)$ of interconnection
(\ref{eq:interconnection}) converges into the set $\Omega_a$
specified by (\ref{eq:domain_of_convergence}). If, in addition,
conditions C1)--C3) of Corollary
\ref{cor:non_uniform_small_gain_attracting_set} hold then the
domain of convergence is given by
(\ref{eq:new_domain_of_convergence}).
\end{cor}

A practically important consequence of this corollary concerns
systems $\mathcal{S}_a$ which are exponentially stable:
\begin{equation}\label{eq:exponentially_stable_dynamics}
\norms{\bfx(t)}\leq \norms{\bfx(t_0)}D_\beta \exp(-\lambda
t)+c\cdot\|h(\bfz(t,\bfz_0))\|_{\infty,[t_0,t]}, \lambda>0, \
D_\beta\geq 1
\end{equation}
In this case the domain (\ref{eq:domain_of_convergence_separable})
of initial conditions ensuring convergence into $\Omega_a$ is
defined as
\[
 D_{\gamma,0}\leq\max_{\kappa>1, \ d\in(0,1)}
-\lambda\left(\ln
\frac{d}{\kappa}\right)^{-1}\frac{\kappa-1}{\kappa}
\frac{h(\bfz_0)}{D_\beta \norms{\bfx_0}+D_\beta \cdot c\cdot
|h(\bfz_0)|\left(1+\frac{\kappa}{1-d}\right)+c|h(\bfz_0)|}
\]

\section{Discussion}\label{sec:discussion}

In this section we discuss some practically relevant outcomes of
the results of Theorem \ref{theorem:non_uniform_small_gain} and
Corollaries \ref{cor:non_uniform_small_gain_attracting_set},
\ref{cor:non_uniform_small_gain_GAS} and their potential
applications to problems of analysis of asymptotic behavior in
nonlinear dynamic systems.

First, in Subsection \ref{subsection:relation_small_gain} we
specify conditions for existence of a trapping region of nonzero
volume in $\Real^n\oplus\Real^m$ in terms of the parameters of
system (\ref{eq:interconnection}) without invoking dependence on
$\bfx(t_0)$, $\bfz(t_0)$ as was done in Theorem
\ref{theorem:non_uniform_small_gain}. The resulting criterion has
a form similar to the standard small-gain conditions
\cite{Zames:66}. The differences and similarities between this new
result and standard small-gain theorems are illustrated with an
example.

Second, in Subsection \ref{subsec:new_adaptive_observers} we
demonstrate how the results of our present contribution can be
applied to address the problem of output nonlinear identification
for systems which cannot be transformed into a canonic observer
form or/and with nonlinear parametrization.

\subsection{Relation to conventional small-gain
theorems}\label{subsection:relation_small_gain}

Conditions specifying state boundedness formulated in Theorem
\ref{theorem:non_uniform_small_gain} and Corollaries
\ref{cor:non_uniform_small_gain_attracting_set},
\ref{cor:non_uniform_small_gain_GAS} depend explicitly on initial
conditions $\bfx(t_0)$, $\bfz(t_0)$. Such dependence is inevitable
when the convergence is allowed to be non-uniform. But if mere
existence of a trapping region is asked for, dependence on initial
conditions may be removed from the statements of the results. The
next corollary presents such modified conditions.

\begin{cor}\label{cor:small_gain_like} Consider
interconnection (\ref{eq:interconnection}) where the system
$\mathcal{S}_a$ satisfies inequality
(\ref{eq:uniformly_attractive_dynamics}) and the function
$\gamma_0(\cdot)$ obeys (\ref{eq:Lipschitz_gamma}). Then there
exists a set $\Omega_\gamma$ of initial conditions corresponding
to the trajectories converging to $\Omega_a$ if the following
condition is satisfied
\begin{equation}\label{eq:small_gain_like_condition}
D_{\gamma,0}\cdot c \cdot \mathcal{G}<1,
\end{equation}
where
\[
\mathcal{G}=
\beta_t^{-1}\left(\frac{d}{\kappa}\right)\frac{k}{k-1}\left(\beta_t(0)\left(1+\frac{\kappa}{1-d}\right)+1\right)
\]
for some $d\in(0,1)$, $\kappa\in(1,\infty)$. In particular,
$\Omega_\gamma$ contains the following domain
\[
\begin{split}
\norms{\bfx(t_0)}\leq \frac{1}{\beta_t(0)}
\left[\frac{1}{D_{\gamma,0}}\left(\beta_t^{-1}\left(\frac{d}{\kappa}\right)\right)^{-1}\frac{k-1}{k}-
c
\left(\beta_t(0)\left(1+\frac{\kappa}{1-d}\right)+1\right)\right]h(\bfz(t_0)).
\end{split}
\]
In case the function $h(\bfz)$ in (\ref{eq:interconnection}) is
continuous, the volume of the set $\Omega_\gamma$ is nonzero in
$\Real^n\oplus\Real^m$.
\end{cor}

Notice that in case the dynamics of the contracting subsystem
$\mathcal{S}_a$ is exponentially stable, i.e. it satisfies
inequality (\ref{eq:exponentially_stable_dynamics}), the term
$\mathcal{G}$ in condition (\ref{eq:small_gain_like_condition})
reduces to
\begin{equation}\label{eq:small_gain_like_condition_exponential}
\mathcal{G}=\frac{1}{\lambda}\cdot
\ln\left(\frac{\kappa}{d}\right)\frac{k}{k-1}\left(D_\beta\left(1+\frac{\kappa}{1-d}\right)+1\right)
\end{equation}
For $D_\beta=1$ the minimal value of $\mathcal{G}$ in
(\ref{eq:small_gain_like_condition_exponential}) can be estimated
as
\begin{equation}\label{eq:small_gain_like_condition_exponential_bound}
\mathcal{G}^\ast=\frac{1}{\lambda} \cdot \min_{d\in(0,1), \
\kappa\in(1,\infty)} \ln\left(\frac{\kappa
}{d}\right)\frac{k}{k-1}\left(2+\frac{\kappa}{1-d}\right)\approx
\frac{15.6886}{\lambda}<\frac{16}{\lambda},
\end{equation}
which leads to an even more simple formulation of
(\ref{eq:small_gain_like_condition_exponential}):
\[
D_{\gamma,0}\cdot \frac{c}{\lambda} \leq \frac{1}{16}
\]

Corollary \ref{cor:small_gain_like} provides an explicit and
easy-to-check condition for existence of a trapping region in the
state space of a class of Lyapunov unstable systems. In addition,
it allows to specify explicitly points $\bfx(t_0)$, $\bfz(t_0)$
which belong to the emergent trapping region. Notice also that the
existence condition, inequality
(\ref{eq:small_gain_like_condition}), has the flavor of
conventional small-gain constraints. Yet, it is substantially
different from these classical results. This is because the
input-output gain for the wandering subsystem, $\mathcal{S}_w$,
may not be finite or need not even be defined.

To elucidate these differences as well as the similarities between
conditions of conventional small-gain theorems and those
formulated in Corollary \ref{cor:small_gain_like} we provide an
example. Consider the following systems
\begin{subequations}
\begin{equation}
\left\{\begin{aligned}
\dot{x}_1&=-\lambda_1 x_1 + c_1 x_2\\
\dot{x}_2&=-\lambda_2 x_2 - c_2 |x_1|
\end{aligned}\right.\label{eq:small_gain_example}
\end{equation}
\begin{equation}
\left\{\begin{aligned}
\dot{x}_1&=-\lambda_1 x_1 + c_1 x_2\\
\dot{x}_2&= - c_2 |x_1|
\end{aligned}\right.\label{eq:non_uniform_small_gain_example}
\end{equation}
\end{subequations}
System (\ref{eq:small_gain_example}) can be viewed as an
interconnection of two input-to-state stable systems, $x_1$ and
$x_2$, with input-output $L_\infty$-gains ${c_1}/{\lambda_1}$ and
$c_2/\lambda_2$ respectively. Therefore, in order to prove state
boundedness of (\ref{eq:small_gain_example}) we can, in principle,
invoke the conventional small-gain theorem. The small-gain
condition in this case is as follows:
\begin{subequations}
\begin{equation}
\frac{c_1}{\lambda_1}\cdot \frac{c_2}{\lambda_2}<1
\end{equation}
The theorem, however, does not apply to system
(\ref{eq:non_uniform_small_gain_example}) because the input-output
gain of its second subsystem, $x_2$, is infinite. Yet, by invoking
Corollary \ref{cor:small_gain_like} it is still possible to show
existence of a weak attracting set in the state space of system
(\ref{eq:non_uniform_small_gain_example}) and specify its basin of
attraction. As follows from Corollary \ref{cor:small_gain_like},
condition
\begin{equation}
\frac{c_1}{\lambda_1} \cdot \frac{c_2}{\lambda_1}<\frac{1}{16}
\end{equation}
ensures existence of the trapping region, and the trapping region
itself is given by
\[
|x_1(t_0)|\leq \left[\frac{1}{c_2}\lambda_1\left(\ln\frac{\kappa
}{d}\right)^{-1}\frac{k-1}{k}-\frac{c_1}{\lambda_1}\left(2+\frac{\kappa}{1-d}\right)\right]
x_2(t_0).
\]
\end{subequations}

\subsection{Output nonlinear identification
problem}\label{subsec:new_adaptive_observers}

In the literature on adaptive control, observation, and
identification a few classes of systems are referred to as {\it
canonic forms} because they guarantee existence of a solution to
the problem and because a large variety of physical models can be
transformed into this class. Among these, perhaps the most widely
known is the {\it adaptive observer canonical form}
\cite{Bastin88}. Necessary and sufficient conditions for
transformation of the original system into this canonical form can
be found, for example, in \cite{Marino90}. These conditions,
however, include restrictive requirements of linearization of
uncertainty-independent dynamics by output injection, and they
also require linear parametrization of the uncertainty.
Alternative approaches \cite{Besancon_2000} heavily rely on
knowledge of the proper Lyapunov function for the
uncertainty-independent part and still assume linear
parametrization.

We now demonstrate how these restrictions can be lifted by
application of our result to the problem of state and parameter
observation. Let us consider systems which can be transformed by
means of static or dynamic feedback\footnote{Notice that
conventional observers in control theory could be viewed as
dynamic feedbacks.} into the following form:
\begin{equation}\label{eq:new_canonical_identifier_form}
\dot{\bfx}=\bff_0(\bfx,t)+\bff(\xivec(t),\thetavec)-\bff(\xivec(t),\hat{\thetavec})+\varevec(t),
\end{equation}
where
\[
 \varevec(t)\in L_{\infty}^m[t_0,\infty], \
 \|\varevec(\tau)\|_{\infty,[t_0,t]}\leq \Delta_\varepsilon
\]
is an external perturbation with known $\Delta_\varepsilon$, and
$\bfx\in\Real^n$. The function
$\xivec:\Real_+\rightarrow\Real^\xi$ is a function of time, which
possibly includes available measurements of the state, and
$\thetavec,\hat{\thetavec}\in\Omega_\theta\subset\Real^d$ are the
unknown and estimated parameters of function the $\bff(\cdot)$,
respectively, and the set $\Omega_\theta$ is bounded.  We assume
that the function $\bff(\xivec(t),\thetavec)$ is locally bounded
in $\thetavec$ uniformly in $\xivec$:
\[
\|\bff(\xivec(t),\thetavec)-\bff(\xivec(t),\hat{\thetavec})\|\leq
D_f\|\thetavec-\hat{\thetavec}\|+\Delta_f
\]
and the values of $D_f\in\Real_+$, $\Delta_f$ are available. The
function $\bff_0(\cdot)$  in
(\ref{eq:new_canonical_identifier_form}) is assumed to satisfy the
following condition.
\begin{assume}\label{assume:observer_error_unperturbed} The system
\begin{equation}\label{eq:observer_unperturbed}
\dot{\bfx}=\bff_0(\bfx,t)+\bfu(t)
\end{equation}
is forward-complete. Furthermore, for all $\bfu(t)$ such that
\[
\|\bfu(t)\|_{\infty,[t_0,t]}\leq
\Delta_u+\|\bfu_0(\tau)\|_{\infty,[t_0,t]}, \ \Delta_u\in\Real_+
\]
there exists a bounded set $\mathcal{A}$, $c>0$ and a function
$\Delta:\Real_+\rightarrow\Real_+$ satisfying the following
inequality
\[
\normss{\bfx(t)}{\Delta(\Delta_u)}\leq
\beta(t-t_0)\normss{\bfx(t_0)}{\Delta(\Delta_u)}+c\|\bfu_0(\tau)\|_{\infty,[t_0,t]}
\]
where $\beta(\cdot):\Real_+\rightarrow\Real_+$,
$\lim_{t\rightarrow\infty}\beta(t)=0$ is a strictly decreasing
function
\end{assume}
Consider the following auxiliary system
\begin{equation}\label{eq:exo_searching}
\dot{\lambdavec}=S(\lambdavec), \
\lambdavec(t_0)=\lambdavec_0\in\Omega_\lambda\subset\Real^\lambda
\end{equation}
where $\Omega_\lambda\subset\Real^n$ is bounded and $S(\lambda)$
is locally Lipschitz. Furthermore, suppose that the following
assumption holds for system (\ref{eq:exo_searching}).
\begin{assume}\label{assume:exo_searching} System (\ref{eq:exo_searching}) is
 Poisson stable in
$\Omega_\lambda$ that is
\[
\forall \ \lambdavec'\in\Omega_\lambda,  \ t'\in\Real_+ \
\Rightarrow \exists t''>t: \
\|\lambdavec(t'',\lambdavec')-\lambdavec'\|\leq \epsilon,
\]
where $\epsilon$ is an arbitrary small positive constant.
Moreover, the trajectory $\lambda(t,\lambdavec_0)$ is dense in
$\Omega_\lambda$:
\[
\forall \lambdavec'\in\Omega_\lambda, \
\epsilon\in\Real_{>0}\Rightarrow \exists \ t\in\Real_+: \
\|\lambdavec'-\lambdavec(t,\lambdavec_0)\|<\epsilon
\]
\end{assume}
Now we are ready to formulate the following statement

\begin{cor}\label{cor:identifier} Consider system
(\ref{eq:new_canonical_identifier_form}) and suppose that the
following conditions hold

C4) the vector-field $\bff_0(\bfx,t)$ in
(\ref{eq:new_canonical_identifier_form}) satisfies Assumption
\ref{assume:observer_error_unperturbed};

C5) there exists a (known) system (\ref{eq:exo_searching})
satisfying Assumption \ref{assume:exo_searching};

C6) there exists  a locally Lipschitz
$\etavec:\Real^\lambda\rightarrow\Real^d$:
\[
\|\etavec(\lambdavec')-\etavec(\lambdavec'')\|\leq D_\eta
\|\lambdavec'-\lambdavec''\|
\]
such that the set $\etavec(\Omega_\lambda)$ is dense in
$\Omega_\theta$;

C7)  system (\ref{eq:new_canonical_identifier_form}) has
steady-state characteristic with respect to the norm
\[
\normss{\cdot}{\Delta(M)}, \ M=2\Delta_f + \Delta_\varepsilon +
\delta
\]
and input $\hat{\thetavec}$, where $\delta$ is  some positive
(arbitrary small) constant.

Consider the following interconnection of
(\ref{eq:new_canonical_identifier_form}),
(\ref{eq:exo_searching}):
\begin{equation}\label{eq:identifier_closed}
\begin{split}
\dot{\bfx}&=\bff_0(\bfx,t)+\bff(\xivec(t),\thetavec)-\bff(\xivec(t),\hat{\thetavec})+\varevec(t)\\
\hat{\thetavec}&=\etavec(\lambdavec)\\
\dot{\lambdavec}&=\gamma \normss{\bfx(t)}{\Delta(M)}S(\lambdavec),
\end{split}
\end{equation}
where $\gamma>0$ satisfies the following inequality
\begin{equation}\label{eq:gamma_identifier}
\begin{split}
\gamma&\leq
\left(\beta_t^{-1}\left(\frac{d}{\kappa}\right)\right)^{-1}\frac{\kappa-1}{\kappa}
\frac{1}{D_\lambda\left(\beta_t(0)
\left(1+\frac{\kappa}{1-d}\right)+1\right)}\\
D_\lambda&=c\cdot D_f\cdot
D_\eta\cdot\max_{\lambdavec\in\Omega\lambda}\|S(\lambdavec)\|
\end{split}
\end{equation}
for some $d\in(0,1)$, $\kappa\in(1,\infty)$. Then,  for
$\lambdavec(t_0)=\lambdavec_0$, some $\thetavec'\in\Omega_\theta$
and all $\bfx(t_0)=\bfx_0\in\Real^n$ the following holds
\begin{equation}\label{eq:identifier_convergence}
\begin{split}
& \lim_{t\rightarrow\infty}\normss{\bfx(t)}{\Delta(M)}=0, \
\lim_{t\rightarrow\infty}\hat{\thetavec}(t)=\thetavec'\in\Omega_\theta
\end{split}
\end{equation}
\end{cor}
Notice that, as has been pointed out in the previous section, in
case the dynamics of (\ref{eq:observer_unperturbed}) is
exponentially stable with rate of convergence equal to $\rho$ and
$\beta(0)=D_\beta$, condition (\ref{eq:gamma_identifier}) will
have the following form
\[
\gamma\leq -\rho\left(\ln
\frac{d}{\kappa}\right)^{-1}\frac{\kappa-1}{\kappa}\frac{1}{D_\lambda\left(D_\beta
\left(1+\frac{\kappa}{1-d}\right)+1\right)}
\]

According to Corollary \ref{cor:identifier}, for the rather
general class of systems (\ref{eq:new_canonical_identifier_form})
it is possible to design an estimator $\hat{\thetavec}(t)$  which
guarantees that not only the "error" vector $\bfx(t)$ reaches a
neighborhood of the origin, but also that the estimates
$\hat{\thetavec}(t)$ converge to some $\thetavec'$ in
$\Omega_\theta$. Both these facts, together with additional
nonlinear persistent excitation conditions
\cite{Cao_2003},\cite{Tyukin:2005:arx}
\[
\exists T>0, \rho\in\mathcal{K}: \ \forall \ \mathcal{T}=[t,t+T],
\ t\in\Real_+ \Rightarrow \exists \tau\in\mathcal{T}:
 |\bff(\xivec(\tau),\thetavec)-\bff(\xivec(\tau),\thetavec')|\geq
\rho(\|\thetavec-{\thetavec}'\|),
\]
in principle allow us to estimate the domain of convergence for
$\hat{\thetavec}(t)$.

Concluding this section we mention that statements of Theorem
\ref{theorem:non_uniform_small_gain} and Corollaries
\ref{cor:non_uniform_small_gain_attracting_set}--\ref{cor:identifier}
constitute additional theoretical tools for the analysis of
asymptotic behavior of systems in cascaded form. In particular
they are complementary to the results of \cite{SIAM:Arcak:2002}
where {\it asymptotic stability} of the following type of systems
\begin{equation}\label{eq:former_iISS}
\begin{split}
\dot{\bfx}&=\bff(\bfx), \\
\dot{\bfz}&=\bfq(\bfx,\bfz),\ \bff:\Real^n\rightarrow\Real^n, \
\bfq:\Real^n\times\Real^m\rightarrow\Real^m\nonumber
\end{split}
\end{equation}
was considered under assumption that the $\bfx$-subsystem is
globally {asymptotically stable} and the $\bfz$-subsystem is
integral input-to-state stable. In contrast to this, our results
apply to establishing {\it asymptotic convergence} for systems
with the following structure
\begin{equation}\label{eq:ours_iISS}
\begin{split}
\dot{\bfx}&=\bff(\bfx,\bfz), \\
\dot{\bfz}&=\bfq(\bfx,\bfz),\
\bff:\Real^n\times\Real^m\rightarrow\Real^n \nonumber
\end{split}
\end{equation}
where the $\bfx$-subsystem is input-to-state stable, and the
$\bfz$-subsystem could be practically integral input-to-state
stable (see Corollary
\ref{cor:non_uniform_small_gain_attracting_set}), although in
general no stability assumptions are imposed on it.



\section{Examples}

In this section we provide two examples of parameter
identification in nonlinearly parameterized systems that cannot be
transformed into the canonical adaptive observer form.

The first example is merely an academical illustration of
Corollary \ref{cor:identifier}, where only one parameter is
unknown and the system itself is a first-order differential
equation. The second example illustrates a possible application of
our results to the problem of identifying the dynamics in living
cells.

{\it Example 1.} Consider the following system
\begin{equation}\label{eq:sin_unknowns}
\dot{x}=-k x + \sin(x\theta +\theta)+u, \ k>0, \ \theta\in[-a,a]
\end{equation}
where $\theta$ is an unknown parameter and $u$ is the control
input. Without loss of generality we let $a=1$, $k=1$. The problem
is to estimate the parameter $\theta$ from measurements of $x$ and
steer the system to the origin. Clearly, the choice
$u=-\sin(x\hat{\theta} +\hat{\theta})$ transforms
(\ref{eq:sin_unknowns}) into
\begin{equation}\label{eq:sin_unknowns_controlled}
\dot{x}=-k x + \sin(x\theta +\theta)-\sin(x\hat{\theta}
+\hat{\theta})
\end{equation}
which satisfies Assumption
\ref{assume:observer_error_unperturbed}. Moreover, the system
\begin{equation}
\begin{split}
\dot{\lambda}_1&=\lambda_1\\
\dot{\lambda}_2&=-\lambda_2, \ \lambda_1^2(t_0)+\lambda_2^2(t_0)=1
\end{split}\nonumber
\end{equation}
with mapping $\etavec=(1, \ 0)^T\lambdavec$ satisfies Assumption
\ref{assume:exo_searching} and therefore
\begin{equation}\label{eq:sin_exo}
\begin{split}
\dot{\lambda}_1&=\gamma |x|\lambda_1\\
\dot{\lambda}_2&=-\gamma|x|\lambda_2, \
\lambda_1^2(t_0)+\lambda_2^2(t_0)=1
\end{split}
\end{equation}
would be a candidate for the control and parameter estimation
algorithm. According to Corollary \ref{cor:identifier}, the goal
will be reached if the parameter $\gamma$ in (\ref{eq:sin_exo})
obeys the following constraint
\[
\gamma\leq -\rho\left(\ln
\frac{d}{\kappa}\right)^{-1}\frac{\kappa-1}{\kappa}\frac{1}{D_\lambda\left(D_\beta
\left(1+\frac{\kappa}{1-d}\right)+1\right)}, \ \rho=k=1, \
D_\beta=1, \ D_\lambda=1
\]
for some $d\in(0,1)$, $\kappa\in(1,\infty)$. Hence, choosing, for
example, $d=0.5$, $\kappa=2$ we obtain that choice
\[
0<\gamma<-\ln\left(\frac{0.5}{2}\right)^{-1}\frac{1}{2}\cdot\frac{1}{6}=0.0601
\]
suffices to ensure that
\[
\lim_{t\rightarrow\infty}x(t)=0, \
\lim_{t\rightarrow\infty}\hat{\theta}(t)=\theta
\]

\begin{figure}
\includegraphics[width=220pt]{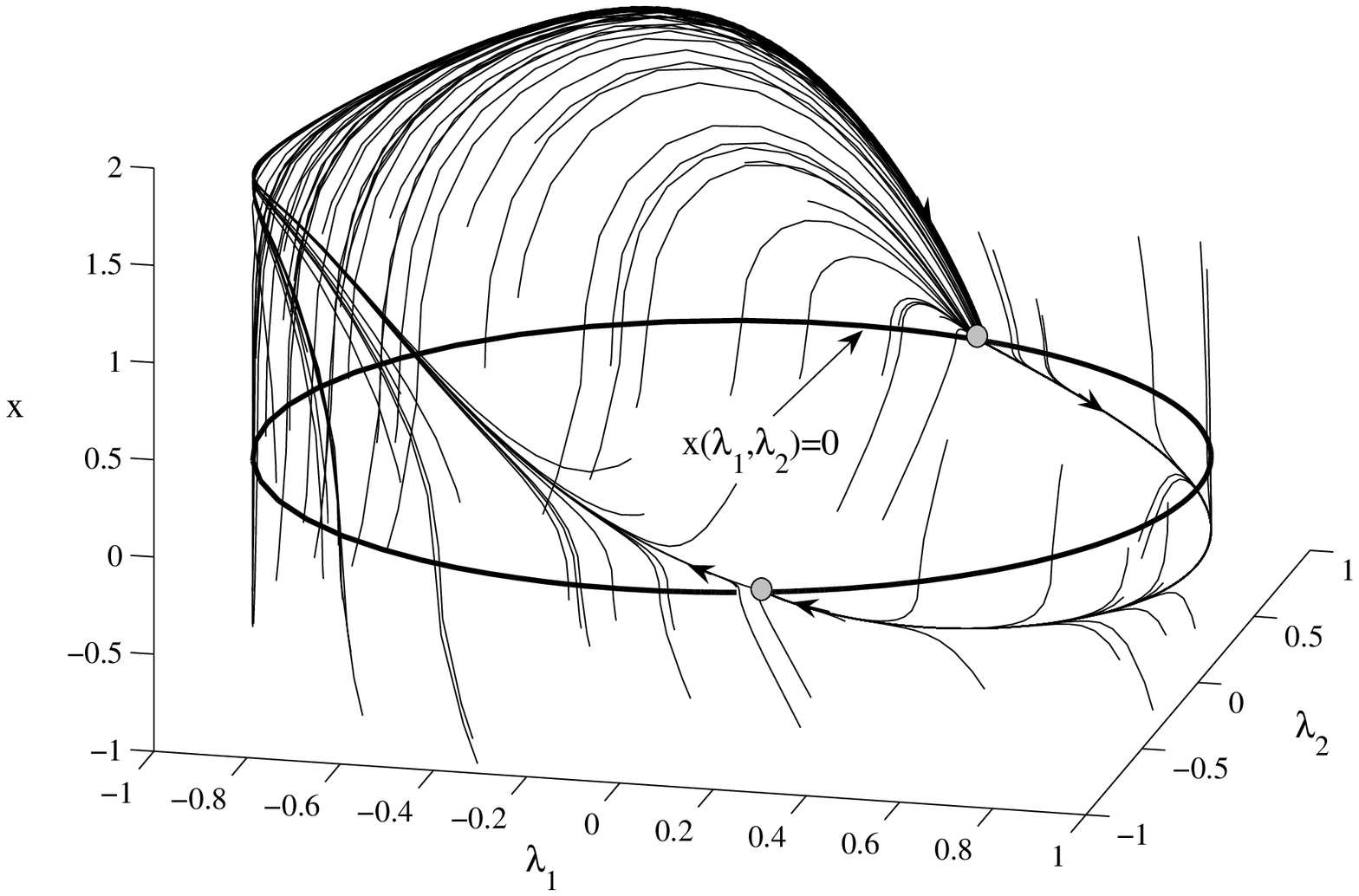}    \hspace{0.56cm}
\includegraphics[width=220pt]{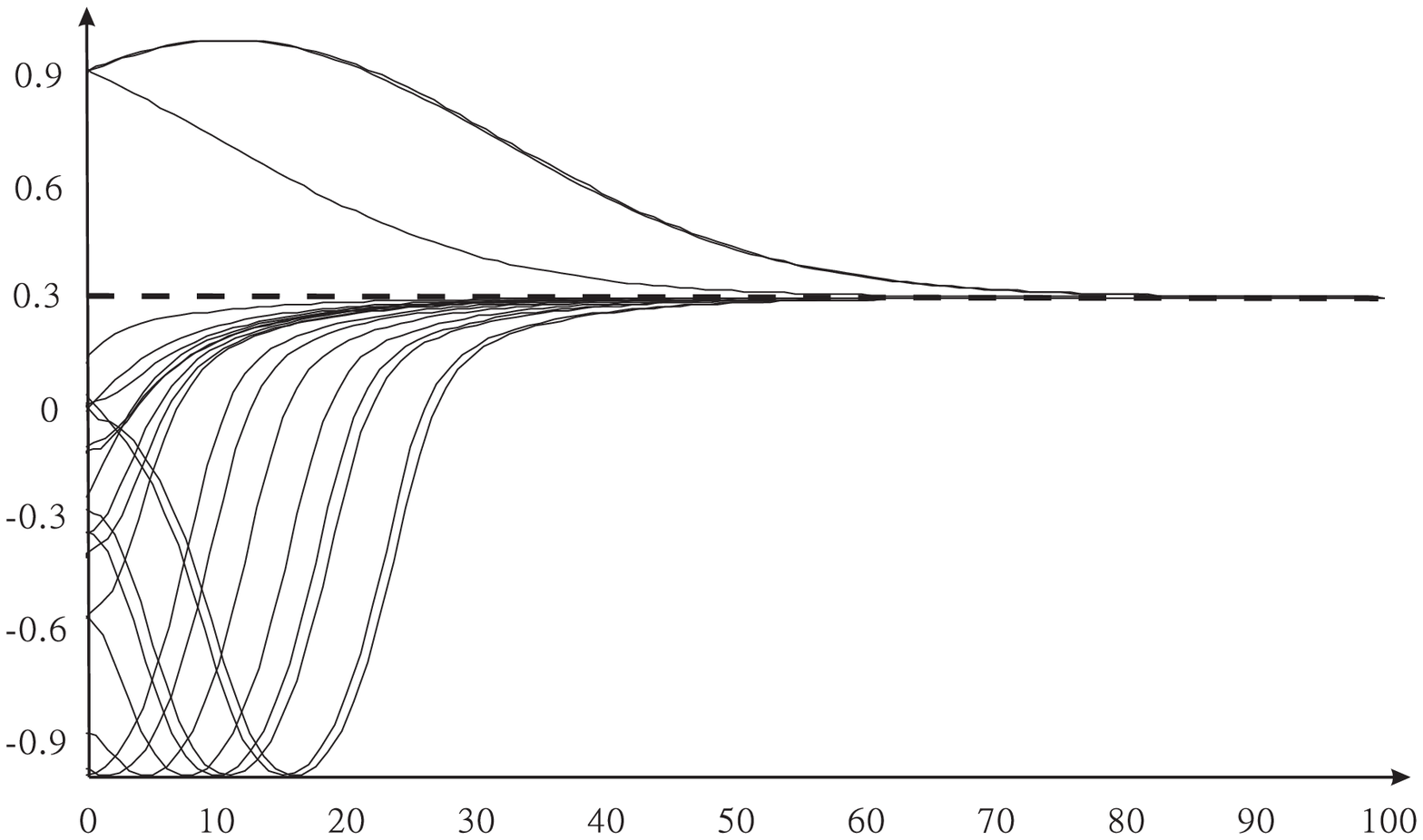}
\begin{center}
\caption{Trajectories of system
(\ref{eq:sin_unknowns_controlled}), (\ref{eq:sin_exo}) (left
panel) and the family of estimates $\hat{\theta}(t)$ of parameter
$\theta$ as functions of time $t$ (right
panel)}\label{fig:adapt_attract}
\end{center}
\end{figure}

We simulated system (\ref{eq:sin_unknowns_controlled}),
(\ref{eq:sin_exo}) with $\theta=0.3$, $\gamma=0.05$ and initial
conditions $x(t_0)$ randomly distributed in  the interval
$[-1,1]$. Results of the simulation are illustrated with Figure
\ref{fig:adapt_attract}, where the phase plots of system
(\ref{eq:sin_unknowns_controlled}), (\ref{eq:sin_exo}) as well as
the trajectories of $\hat{\theta}(t)$ are given.

{\it Example 2.} Consider the problem of modelling electrical
activity in biological cells from the input-output data in current
clamp experiments.  The simplest mathematical model, which
captures a fairly large variety of phenomena like periodic
bursting in response to constant stimulation is the classical
Hindmarsh and Rose model neuron without adaptation currents
\cite{Nature:HR:1982}:
\begin{equation}\label{eq:HR}
 \begin{split}
\dot{x}_1 &= -a x_1^3 + b x_1^2 + x_2  + \alpha u\\
\dot{x}_2 &= c - \beta x_2 - d x_1^2
\end{split}
\end{equation}
where variable $x_1$ is the membrane potential, $x_2$ stands for
the ionic currents in the cell, $u$ is the input current, and $a$,
$b$, $c$, $d$, $\alpha$, $\beta$ $\in\Real$ are parameters. While
the parameters of the first equation can, in principle, be
identified experimentally by blocking the ionic channels in the
cells and measuring the membrane conductance, identification of
parameters $\beta$, $d$ is a difficult problem, as information
about ionic currents $x_2$ is rarely available.

Conventional techniques \cite{Bastin88} cannot be applied directly
to this problem as the model (\ref{eq:HR}) is not in canonical
adaptive observer form. Let us illustrate how our results can be
used to derive the unknown parameters of (\ref{eq:HR}) such that
the reconstructed model fits the observed data. Assume, first,
that parameters $a$, $b$, $c$, $\alpha$ in the first equation of
(\ref{eq:HR}) are known, whereas parameters $\beta$, $d$ in the
second equation are unknown. This corresponds to the realistic
case where the time constant of current $x_2$ and coupling between
$x_1$ and $x_2$ are uncertain. In our example we assumed that
\[
\beta\in\Omega_\beta=[0.3,0.7], \ \ d\in\Omega_d=[2,3], \ a=1, \
b=3, \ \alpha=0.7, \ c=0.5
\]

As a candidate for the observer  we select the following system
\begin{equation}\label{eq:observer}
\begin{split}
\dot{\hat{x}}&=\rho(x_1- \hat{x}) -a x_1^3 + b x_1^2+ \alpha u +
f(\hat{\beta},\hat{d},t), \ \rho\in\Real_{>0}\\
\end{split}
\end{equation}
where $\hat{\beta}$, $\hat{d}$ are parameters to be adjusted and
the function $f(\hat{\beta},\hat{d},t)$ is specified as
\[
f(\hat{\beta},\hat{d},t)=\int_{0}^t
e^{-\hat{\beta}(t-\tau)}(\hat{d}x_1^2(\tau)+c)d\tau
\]
Then the dynamics of $\tilde{x}(t)=x(t)-\hat{x}(t)$ satisfies the
following differential equation
\[
\dot{\tilde{x}}=-\rho \tilde{x}+
f(\beta,d,t)-f(\hat{\beta},\hat{d},t)
\]

The function $f(\beta,d,t)$ satisfies the following inequality
\begin{equation}
\begin{split}
|f(\beta,d,t)-f(\hat{\beta},\hat{d},t)|&\leq
|f(\beta,d,t)-f(\hat{\beta},{d},t)|+|f(\hat{\beta},d,t)-f(\hat{\beta},\hat{d},t)|\\
&\leq D_{f,\beta}|\beta-\hat{\beta}|+D_{f,d}|d-\hat{d}|+
\epsilon(t),
\end{split}\nonumber
\end{equation}
where $\epsilon(t)$ is an exponentially decaying term, and
\begin{equation}\label{eq:example_bounds_Df}
D_{f,\beta}=\max_{\hat{\beta},\beta\in\Omega_\beta, \
d\in\Omega_d}\left\{\frac{1}{\beta\hat{\beta}}(d\|x_1(\tau)\|_{\infty,[t_0,\infty]}+c)\right\},
\
D_{f,d}=\max_{\hat{\beta}\in\Omega_\beta}\left\{\frac{1}{\hat{\beta}}\|x_1(\tau)\|_{\infty,[t_0,\infty]}\right\}
\end{equation}
Furthermore, Assumption \ref{assume:observer_error_unperturbed} is
satisfied for system
\begin{equation}\label{eq:example:assumption3}
\dot{\tilde{x}}=-\rho\tilde x + \upsilon(t),
\end{equation}
with
\[
\Delta(\Delta_u)=\frac{\Delta_u}{\rho}.
\]
In particular, for all  $\upsilon(t): \
\|\upsilon(\tau)\|_{\infty,[t_0,t]}\leq
\Delta_u+\|\upsilon_0(\tau)\|_{\infty,[t_0,t]}$ the following
inequality holds:
\begin{equation}\label{eq:example:assumption3_1}
\|x(t)\|_{\Delta(\Delta_u)}\leq e^{-\rho
(t-t_0)}\|x(t_0)\|_{\Delta(\Delta_u)}+\frac{1}{\rho}\|\upsilon_0(\tau)\|_{\infty,[t_0,t]}.
\end{equation}
To see this consider the general solution of
(\ref{eq:example:assumption3}):
\[
x(t)=e^{-\rho (t-t_0)}x(t_0) + e^{-\rho t}\int_{t_0}^t e^{\rho
\tau}\upsilon(\tau)d\tau
\]
and derive an estimate of $|x(t)|$. This estimate has the
following form:
\[
\begin{split}
|x(t)|&\leq e^{-\rho
(t-t_0)}|x(t_0)|+\frac{1}{\rho}\left(1-e^{-\rho(t-t_0)}\right)
\|\upsilon(\tau)\|_{\infty,[t_0,t]}\\
&\leq e^{-\rho
(t-t_0)}\left(|x(t_0)|-\frac{1}{\rho}\Delta_u\right)+\frac{1}{\rho}
\left(\|\upsilon_0(\tau)\|_{\infty,[t_0,t]}+\Delta_u\right)\\
&\leq e^{-\rho (t-t_0)}\|x(t_0)\|_{\Delta(\Delta_u)}+
\frac{1}{\rho}\left(\|\upsilon_0(\tau)\|_{\infty,[t_0,t]}+\Delta_u\right)
\end{split}
\]
Hence
\[
|x(t)|-\frac{1}{\rho}\Delta_u\leq  e^{-\rho
(t-t_0)}\|x(t_0)\|_{\Delta(\Delta_u)} + \frac{1}{\rho}
\|\upsilon_0(\tau)\|_{\infty,[t_0,t]},
\]
which automatically implies (\ref{eq:example:assumption3_1}).

Let us define subsystem (\ref{eq:exo_searching}). Consider the
following system of differential equations
\begin{equation}\label{eq:example_exo_system}
\begin{split}
\dot{\lambda}_1&=\lambda_2\\
\dot{\lambda}_2&=-\omega_1^2\lambda_1\\
\dot{\lambda}_3&=\lambda_4\\
\dot{\lambda}_4&=-\omega_2^2\lambda_3, \ \lambdavec_0=(1,0,1,0)^T
\end{split}
\end{equation}
where $\Omega_\lambda$ is the $\omega$-limit set of the point
$\lambdavec_0$, and $\omega_1, \ \omega_2\in\Real$.  System
(\ref{eq:example_exo_system}), therefore, satisfies Assumption
\ref{assume:exo_searching}. Given that domains $\Omega_\beta$,
$\Omega_d$ are known,  select
\begin{equation}\label{eq:example_eta_map}
\begin{split}
\etavec:& \ \Real^n\rightarrow\Real^2,  \ \etavec=(\eta_1(\lambdavec),\eta_2(\lambdavec))\\
\hat{\beta}=\eta_1(\lambdavec)&=\frac{1}{2}\left(\frac{2
\arcsin(\lambda_1)}{\pi}+1\right)\cdot 0.4 + 0.3, \
\hat{d}=\eta_2(\lambdavec)=\frac{1}{2}\left(\frac{2\arcsin(\lambda_3)}{\pi}+1\right)+2
\end{split}
\end{equation}
Choosing
\[
\frac{\omega_1}{\omega_2}=\pi
\]
we ensure that $\etavec(\Omega_\lambda)$ is dense in
$\Omega_\beta\times\Omega_d$. Given that $\hat{\beta}$, $\hat{d}$
are bounded and $\hat{\beta}\geq 0.3$, $D_{f,\beta}$ and $D_{f,d}$
in (\ref{eq:example_bounds_Df}) are also bounded because for the
given range of parameters signal $x_1(t)$ is always bounded.
Hence, according to Corollary \ref{cor:identifier},
interconnection of (\ref{eq:observer}), (\ref{eq:example_eta_map})
and
\begin{equation}\label{eq:example_exo_system_1}
\begin{split}
\dot{\lambda}_1&=\gamma \|\tilde{x}(t)\|_{\Delta(\delta)}\cdot\lambda_2\\
\dot{\lambda}_2&=-\gamma \|\tilde{x}(t)\|_{\Delta(\delta)}\cdot\omega_1^2\lambda_1\\
\dot{\lambda}_3&=\gamma \|\tilde{x}(t)\|_{\Delta(\delta)}\cdot\lambda_4\\
\dot{\lambda}_4&=-\gamma
\|\tilde{x}(t)\|_{\Delta(\delta)}\cdot\omega_2^2\lambda_3, \
\lambdavec_0=(1,0,1,0)^T\nonumber
\end{split}
\end{equation}
with arbitrary small $\delta>0$ and properly chosen $\gamma>0$
ensures that
\[
\lim_{t\rightarrow\infty}{\|\tilde{x}(t)\|_{\Delta(\delta)}}=0, \
\lim_{t\rightarrow\infty}\hat{\beta}(t)=\beta'\in\Omega_\beta, \
 \lim_{t\rightarrow\infty}\hat{d}(t)=d'\in\Omega_d
\]
This in turn implies a successful fit of the model to the
observations.

\begin{figure}
\begin{center}
\begin{minipage}[h]{0.45\linewidth}
\includegraphics[width=\textwidth]{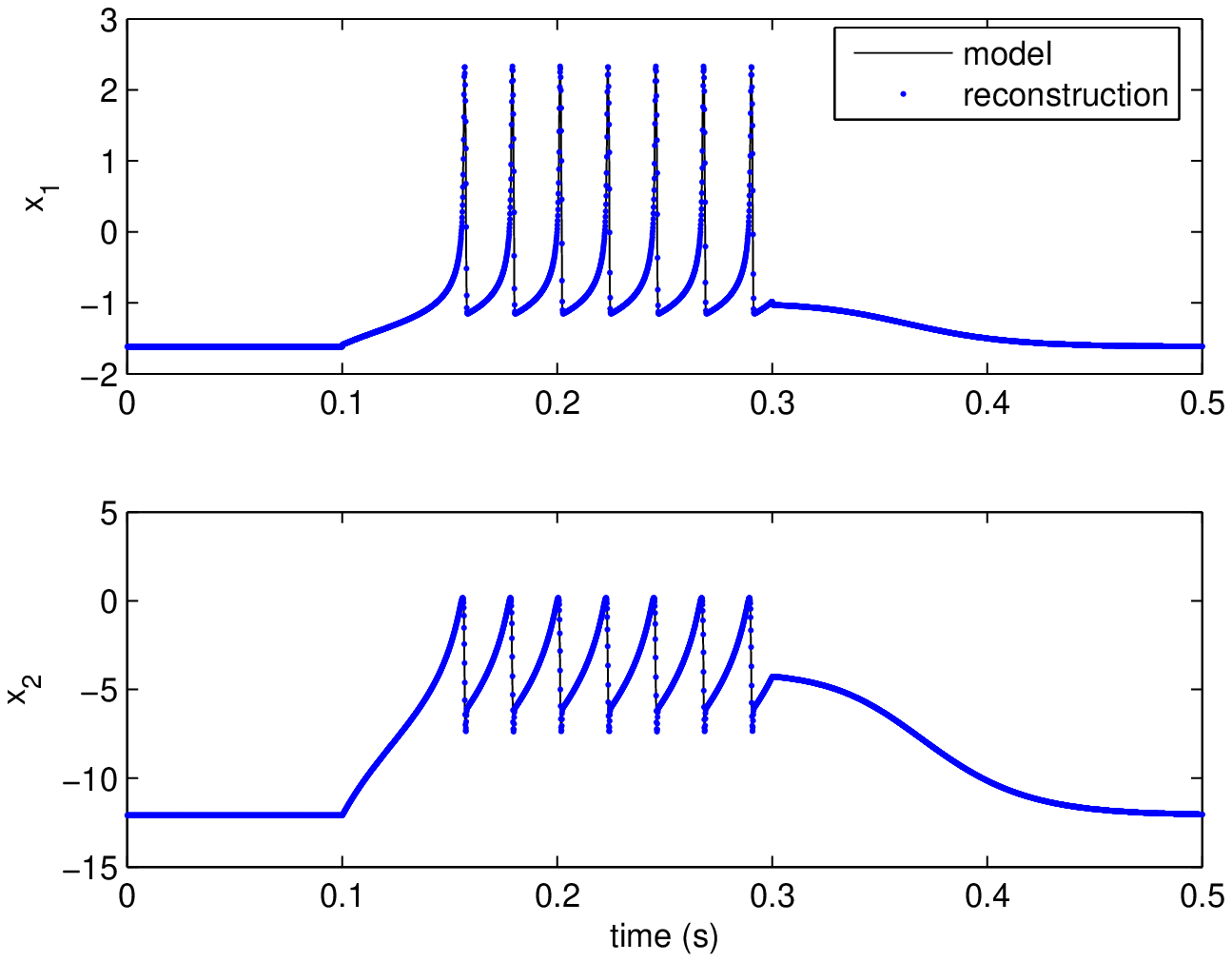}
\end{minipage}
\hspace{7mm}
\begin{minipage}[h]{0.45\linewidth}
\includegraphics[width=\textwidth]{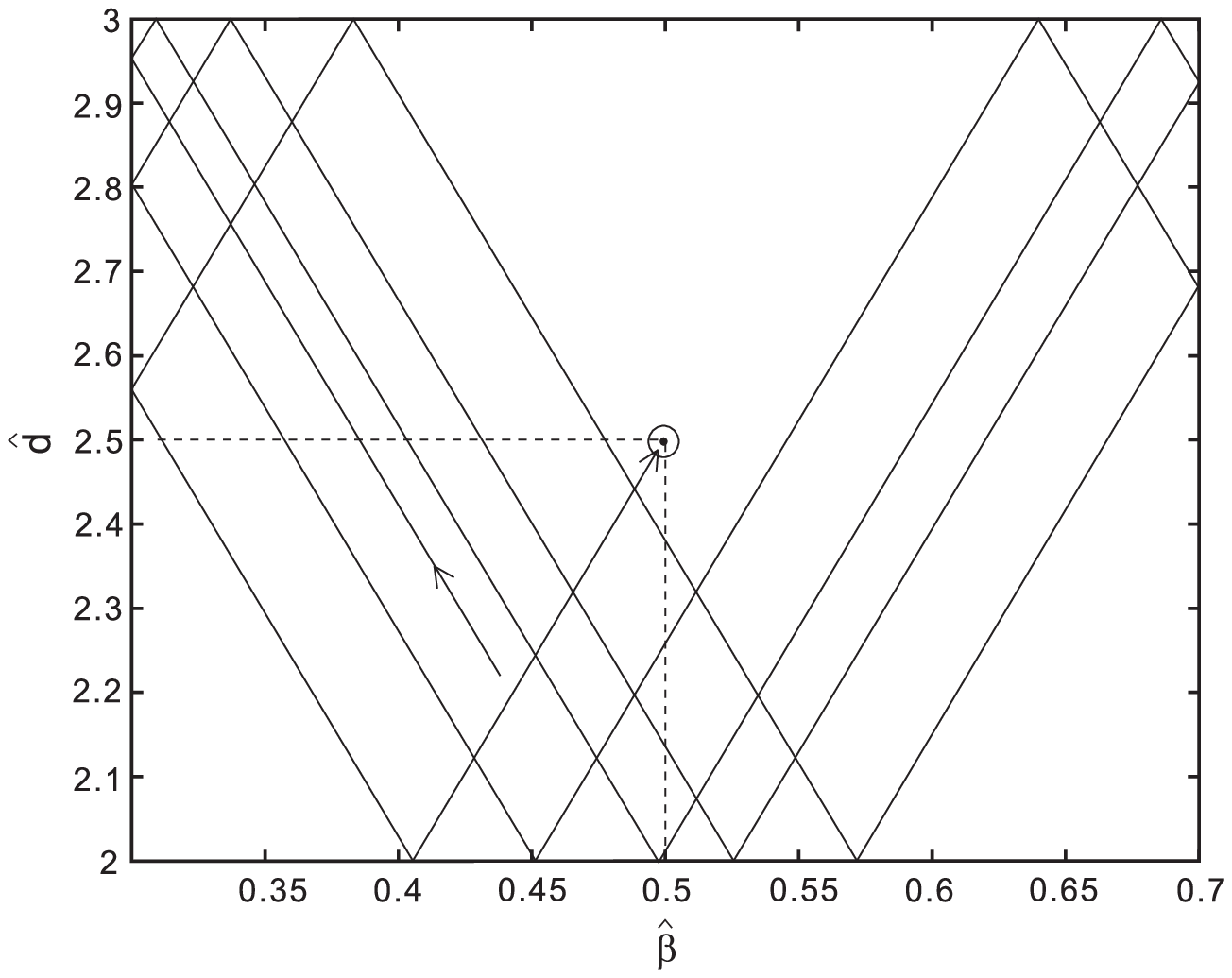}
\end{minipage}
\caption{Left panel -- trajectories $x_1(t)$, $x_2(t)$ of system
(\ref{eq:HR}) plotted for the nominal values of parameters
$\beta=0.5$, $d=2.5$ (model), and for the  values
$\beta=\hat{\beta}(t_0+T)$, $d=\hat{d}(t_0+T)$, where $T$ is the
total simulation time (reconstruction). Input $u(t)$ is a
rectangular impulse with amplitude $0.7$ starting at $t=100$ and
ending at $t=300$. Right panel -- searching dynamics in the
bounded parameter space (a segment of the trajectory
$\hat{\beta}(t),\hat{d}(t)$ towards the end of the
simulation).}\label{fig:trajectories}
\end{center}
\end{figure}

We simulated the system with $\rho=10$ and $\gamma=3\cdot 10^{-4}$
for $\beta=0.5$, $d=2.5$. The results of the simulations are
provided in  figure \ref{fig:trajectories}. It can be seen from
this figure that the reconstruction is successful and the
parameters converge into a small neighborhood of the actual
values.

\section{Conclusion}

We proposed tools for the analysis of asymptotic behavior of a
class of dynamical systems. In particular, we consider an
interconnection of an input-to-state stable system with an
unstable or integrally input-to-state dynamics. Our results allow
to address a variety of problems in which convergence may not be
unform with respect to initial conditions. It is necessary to
notice that the proposed method does not require complete
knowledge of the dynamical systems in question. Only qualitative
information like, for instance, characterization of input-to-state
stability of is necessary for application of our results. We
demonstrated how our analysis can be used in the problems of
synthesis and design -- in particular to  problems of nonlinear
regulation and parameter identification of nonlinear parameterized
systems. The examples show the relevance of our approach in those
domains where application of the standard techniques is either not
possible or too complicated.

\section{Acknowledgment}

The authors are thankful to  Peter Jurica and Tatiana Tyukina for
their enthusiastic help and comments during the preparation of
this manuscript.

\section{Appendix}

{\it Proof of Theorem \ref{theorem:non_uniform_small_gain}.} Let
the conditions of the theorem be satisfied for  given
$t_0\in\Real_+$: $\bfx(t_0)=\bfx_0$, $\bfz(t_0)=\bfz_0$. Notice
that in this case $h(\bfz_0)\geq 0$, otherwise requirement
(\ref{eq:non_uniform_small_gain}) will be violated. Consider the
sequence (\ref{eq:shrinking_volumes}) of volumes $\Omega_i$
induced by $\mathcal{S}$:
\[
\Omega_i=\{\bfx\in\mathcal{X}, \ \bfz\in\mathcal{Z}| \
h(\bfz(t))\in H_i \}
\]
To prove the theorem we show that $0\leq h(\bfz(t))\leq h(\bfz_0)$
for all $t\geq t_0$. For the given partition
(\ref{eq:shrinking_volumes}) we consider two alternatives.

First, in the degenerative case, the state $\bfx(t)\oplus\bfz(t)$
enters some $\Omega_j$, $j\geq 0$ and stays there for all $t\geq
t_0$ which automatically guarantees that $0\leq|h(\bfz)|\leq
h(\bfz_0)$. Then, according to (\ref{eq:attracting}) the
trajectory $\bfx(t)$ satisfies the following inequality:
\begin{equation}\label{eq:theorem_non_uniform_proof_1}
\norms{\bfx(t)}\leq
\beta(\norms{\bfx_0},t-t_0)+c\|h(\bfz(t))\|_{\infty,[t_0,t]} \leq
\beta(\norms{\bfx_0},t-t_0)+c|h(\bfz_0)|
\end{equation}
Taking into account that $\beta(\cdot,\cdot)\in\mathcal{KL}$ we
can conclude that (\ref{eq:theorem_non_uniform_proof_1}) implies
that
\begin{equation}\label{eq:theorem_non_uniform_proof_2}
\lim\sup_{t\rightarrow\infty}\norms{\bfx(t)}=c|h(\bfz_0)|
\end{equation}
Therefore the statements of the theorem hold.

Let us consider the second alternative, where the state
$\bfx(t)\oplus\bfz(t)$ does not belong to $\Omega_j$ for all
$t\geq t_0$. Given that $h(\bfz(t))$ is monotone and
non-increasing in $t$, this implies that there exists an ordered
sequence of time instants $t_j$:
\begin{equation}\label{eq:theorem_non_uniform_proof_3}
t_0>t_1>t_2\cdots t_j>t_{j+1} \cdots
\end{equation}
such that
\begin{equation}\label{eq:theorem_non_uniform_proof_4}
 h(\bfz(t_i))=\sigma_i h(\bfz_0)
\end{equation}
Hence in order to prove the theorem we must show that the sequence
$\{t_i\}_{i=0}^\infty$ does not converge. In other words, the
boundary $\sigma_\infty h(\bfz_0)=0$ will not be reached in finite
time.

In order to do this let us estimate the upper bounds for the
following differences
\[
 T_i=t_{i+1}-t_i
\]
Taking into account inequality (\ref{eq:integral}) and the fact
that $\gamma_0(\cdot)\in\mathcal{K}_e$ we can derive that
\begin{equation}\label{eq:theorem_non_uniform_proof_5}
h(\bfz(t_{i}))-h(\bfz(t_{i+1}))\leq
T_i\max_{\tau\in[t_{i},t_{i+1}]}\gamma_0(\norms{\bfx(\tau)})\leq
T_i\gamma_0(\normsinf{\bfx(\tau)}{[t_i,t_{i+1}]})
\end{equation}
According to the definition of $t_i$ in
(\ref{eq:theorem_non_uniform_proof_4}) and noticing that the
sequence $\mathcal{S}$ is strictly decreasing we have
\[
h(\bfz(t_{i}))-h(\bfz(t_{i+1}))=(\sigma_{i}-\sigma_{i+1})h(\bfz_0)>0
\]
Hence $h(\bfz_0)>0$ implies that
$\gamma_0(\normsinf{\bfx(\tau)}{[t_i,t_{i+1}]})>0$ and, therefore,
(\ref{eq:theorem_non_uniform_proof_5}) results in
 the following estimate of $T_i$
\begin{equation}\label{eq:theorem_non_uniform_proof_6}
T_i\geq
\frac{h(\bfz(t_{i}))-h(\bfz(t_{i+1}))}{\gamma_0(\normsinf{\bfx(\tau)}{[t_i,t_{i+1}]})}=\frac{h(\bfz_0)(\sigma_i-\sigma_{i+1})}{\gamma_0(\normsinf{\bfx(\tau)}{[t_i,t_{i+1}]})}
\end{equation}
Taking into account that $h(\bfz(t))$ is non-increasing over
$[t_i,t_{i+1}]$ and using (\ref{eq:attracting}) we can bound the
norm $\normsinf{\bfx(\tau)}{[t_i,t_{i+1}]}$ as follows
\begin{equation}\label{eq:theorem_non_uniform_proof_7}
\normsinf{\bfx(\tau)}{[t_i,t_{i+1}]}\leq
\beta(\norms{\bfx(t_i)},0)+c\|h(\bfz(\tau))\|_{\infty,[t_i,t_{i+1}]}\leq
\beta(\norms{\bfx(t_i)},0)+c\cdot \sigma_i h(\bfz_0)
\end{equation}
Hence, combining (\ref{eq:theorem_non_uniform_proof_6}) and
(\ref{eq:theorem_non_uniform_proof_7}) we  obtain that
\[
T_i\geq
\frac{h(\bfz_0)(\sigma_i-\sigma_{i+1})}{\gamma_0(\sigma_i(\sigma_i^{-1}\beta(\norms{\bfx(t_i)},0)+c\cdot
h(\bfz_0)))}
\]
Then, using property (\ref{eq:assume:gamma_0}) of function
$\gamma_0$ we can derive that
\begin{equation}\label{eq:theorem_non_uniform_proof_8}
T_i\geq
\frac{h(\bfz_0)(\sigma_i-\sigma_{i+1})}{\gamma_{0,1}(\sigma_i)}\frac{1}{\gamma_{0,2}(\sigma_i^{-1}\beta(\norms{\bfx(t_i)},0)+c\cdot
h(\bfz_0)))}
\end{equation}
Taking into account condition (\ref{eq:tau_system_2}) of the
theorem, the theorem will be proven if we assure that
\begin{equation}\label{eq:theorem_non_uniform_proof_9}
T_i\geq \tau_i
\end{equation}
for all $i=0,1,2,\dots,\infty$. We prove this claim by induction
with respect to the index $i=0,1,\dots,\infty$. We start with
$i=0$, and then show that for all $i> 0$ the following implication
holds
\begin{equation}\label{eq:theorem_non_uniform__proof_9_1}
T_i\geq \tau_i \Rightarrow T_{i+1}\geq \tau_{i+1}
\end{equation}

Let us prove that  (\ref{eq:theorem_non_uniform_proof_9}) holds
for $i=0$. To this purpose consider the term
$({\sigma_{i}-\sigma_{i+1}})/{\gamma_{0,1}(\sigma_i)}$. As follows
immediately from  the conditions of the theorem, equation
(\ref{eq:tau_system_1}), we have that
\begin{equation}\label{eq:theorem_non_uniform_proof_10}
\frac{\sigma_{i}-\sigma_{i+1}}{\gamma_{0,1}(\sigma_i)}\geq \tau_i
\Delta_0 \ \forall \ i\geq 0
\end{equation}
In particular
\[
\frac{\sigma_{0}-\sigma_{1}}{\gamma_{0,1}(\sigma_0)}\geq \tau_0
\Delta_0
\]
Therefore, inequality (\ref{eq:theorem_non_uniform_proof_8})
reduces to
\begin{equation}\label{eq:theorem_non_uniform_proof_11}
T_0\geq \tau_0\Delta_0
\frac{h(\bfz_0)}{\gamma_{0,2}(\sigma_0^{-1}\beta(\norms{\bfx(t_0)},0)+c\cdot
h(\bfz_0))}
\end{equation}
Moreover, taking into account Condition \ref{assume:phi_system}
and (\ref{eq:phi_system}), (\ref{eq:u_system}) we can derive the
following estimate:
\[
\sigma_0^{-1}\beta(\norms{\bfx(t_0)},0)\leq
\sigma_0^{-1}\phi_0(\norms{\bfx(t_0)})+\sigma_0^{-1}\upsilon_0(c\cdot|h(\bfz_0)|\sigma_0)\leq
B_1(\norms{\bfx_0})+B_2(|h(\bfz_0)|,c)
\]
According to the theorem conditions $\bfx_0$ and $\bfz_0$ satisfy
inequality  (\ref{eq:non_uniform_small_gain}). This in turn
implies that
\begin{equation}\label{eq:theorem_non_uniform_proof_12}
\gamma_{0,2}(\sigma_0^{-1}\beta(\norms{\bfx(t_0)},0)+c\cdot
h(\bfz_0))\leq
\gamma_{0,2}(B_1(\norms{\bfx_0})+B_2(|h(\bfz_0)|,c)+c\cdot
h(\bfz_0))\leq \Delta_0\cdot h(\bfz_0)
\end{equation}
Combining   (\ref{eq:theorem_non_uniform_proof_11}) and
(\ref{eq:theorem_non_uniform_proof_12}) we obtain the desired
inequality
\[
T_0\geq \tau_0\Delta_0
\frac{h(\bfz_0)}{\gamma_{0,2}(\sigma_0^{-1}\beta(\norms{\bfx(t_0)},0)+c\cdot
h(\bfz_0))}\geq \tau_0\frac{\Delta_0 h(\bfz_0)}{\Delta_0
h(\bfz_0)}=\tau_0
\]
Thus the basis of induction is proven.

Let us assume that (\ref{eq:theorem_non_uniform_proof_9}) holds
for all $i=0,\dots,n$, $n\geq 0$. We shall prove now that
implication (\ref{eq:theorem_non_uniform__proof_9_1}) holds for
$i=n+1$. Consider the term $\beta(\norms{\bfx(t_{n+1})},0)$:
\begin{equation}
\begin{split}
\beta(\norms{\bfx(t_{n+1})},0)&\leq
\beta(\beta(\norms{\bfx(t_n)},T_n)+c\|h(\bfz(\tau))\|_{\infty,[t_{n},t_{n+1}]},0)\\
& \leq \beta(\beta(\norms{\bfx(t_n)},T_n)+c\cdot \sigma_n \cdot
h(\bfz_0),0)
\end{split}\nonumber
\end{equation}
Taking into account Condition
\ref{assume:desired_contraction_rates_1} (specifically, inequality
(\ref{eq:tau_system})) and
(\ref{eq:phi_system})--(\ref{eq:rho_system}) we can derive that
\begin{equation}\label{eq:theorem_non_uniform_proof_13}
\beta(\norms{\bfx(t_{n+1})},0)\leq
\beta(\xi_n\cdot\beta(\norms{\bfx(t_n)}),0)+c\cdot \sigma_n \cdot
h(\bfz_0),0)\leq \phi_1(\norms{\bfx(t_n)})+\upsilon_1(c\cdot
|h(\bfz)_0|\cdot \sigma_n)
\end{equation}
Notice that, according to the inductive hypothesis ($T_i\geq
\tau_i$), the following holds
\begin{equation}\label{eq:theorem_non_uniform_proof_14}
\norms{\bfx(t_{i+1})}\leq
\beta(\norms{\bfx(t_{i})},T_{i})+c\cdot\sigma_{i}\cdot
h(\bfz_0)\leq
\xi_{i}\beta(\norms{\bfx(t_{i})},0)+c\cdot\sigma_{i}\cdot
h(\bfz_0)
\end{equation}
for all $i=0,\dots,n$. Then
(\ref{eq:theorem_non_uniform_proof_13}),
(\ref{eq:theorem_non_uniform_proof_14}),
(\ref{eq:phi_system})--(\ref{eq:rho_system}) imply that
\begin{equation}\label{eq:theorem_non_uniform_proof_15}
\begin{split}
& \beta(\norms{\bfx(t_{n+1})},0)\leq
\phi_1(\xi_{i}\beta(\norms{\bfx(t_{n-1})},0)+c\cdot\sigma_{n-1}\cdot
h(\bfz_0))+ \upsilon_1(c\cdot |h(\bfz)_0|\cdot \sigma_n)\\
&\leq
\phi_2(\norms{\bfx(t_{n-1})})+\upsilon_2(c\cdot|h(\bfz_0)|\cdot\sigma_{n-1})+\upsilon_1(c\cdot|h(\bfz_0)|\cdot\sigma_{n})\\
&\leq
\phi_{n+1}(\norms{\bfx_0})+\sum_{i=1}^{n+1}\upsilon_{i}(c\cdot|h(\bfz_0)|\sigma_{n+1-i})\leq
\phi_{n+1}(\norms{\bfx_0})+\sum_{i=0}^{n+1}\upsilon_{i}(c\cdot|h(\bfz_0)|\sigma_{n+1-i})
\end{split}
\end{equation}
According to Condition \ref{assume:phi_system}, term
\[
\sigma_{n+1}^{-1}\left(
\phi_{n+1}(\norms{\bfx_0})+\sum_{i=0}^{n+1}\upsilon_{i}(c\cdot|h(\bfz_0)|\sigma_{n+1-i})\right)
\]
is bounded from above by the sum
\[
B_1(\norms{\bfx_0})+B_2(|h(\bfz_0)|,c)
\]
Therefore, monotonicity of $\gamma_{0,2}$, estimate
(\ref{eq:theorem_non_uniform_proof_15}), and inequality
(\ref{eq:non_uniform_small_gain}) lead to the following inequality
\[
\gamma_{0,2}(\sigma_{n+1}^{-1}\beta(\norms{\bfx(t_{n+1}}),0)+c\cdot
h(\bfz_0))\leq
\gamma_{0,2}(B_1(\norms{\bfx_0})+B_2(|h(\bfz_0)|,c)+c\cdot
h(\bfz_0))\leq h(\bfz_0)\Delta_0
\]
Hence, according to (\ref{eq:theorem_non_uniform_proof_8}),
(\ref{eq:theorem_non_uniform_proof_10}) we have:
\[
T_{n+1}\geq
\frac{(\sigma_{n+1}-\sigma_{n+2})}{\gamma_{0,1}(\sigma_{n+1})}\frac{h(\bfz_0)}{\gamma_{0,2}(\sigma_{n+1}^{-1}\beta(\norms{\bfx(t_{n+1})},0)+c\cdot
h(\bfz_0))}\geq \tau_{n+1}\frac{\Delta_0 h(\bfz_0)}{\Delta_0
h(\bfz_0)}=\tau_{n+1}
\]
Thus implication (\ref{eq:theorem_non_uniform__proof_9_1}) is
proven. This implies that $h(\bfz(t))\in[0,h(\bfz_0)]$ for all
$t\geq t_0$ and, consequently, that
(\ref{eq:theorem_non_uniform_proof_2}) holds.  {\it The theorem is
proven.}

\vskip 2mm

{\it Proof of Lemma \ref{lem:steady_state_convergence}.} As
follows from the assumptions, $h(\bfz(t,\bfz_0))$ is bounded.
Assume it belongs to the following interval $[a,h(\bfz_0)]$,
$a\leq h(\bfz_0).$ Therefore, as follows from (\ref{eq:integral})
we can conclude that
\begin{equation}\label{eq:cor:steady_state_1}
0\leq \int_{t_0}^\infty
\gamma_1(\norms{\bfx(\tau,\bfx_0)})d\tau\leq
h(\bfz_0)-h(\bfz(t,\bfz_0))\leq\infty
\end{equation}
On the other hand, taking into account that $h(\bfz(t,\bfz_0))$ is
bounded and monotone in $t$ (every subsequence of which is this is
again monotone) and applying the Bolzano-Weierstrass theorem we
can conclude that $h(\bfz(t,\bfz_0))$ converges in
$[a,h(\bfz_0)]$. In particular, there exists $\bar{h}\in
[a,h(\bfz_0)]$ such that
\begin{equation}\label{eq:cor:steady_state_1_1}
\lim_{t\rightarrow\infty} h(\bfz(t,\bfz_0)) = \bar{h}
\end{equation}
According to the lemma assumptions, system $\mathcal{S}_a$ has
steady-state characteristics. This means that there exists a
constant $\bar{x}\in\Real_+$ such that
\begin{equation}\label{eq:cor:steady_state_2}
\lim_{t\rightarrow\infty}\norms{\bfx(t,\bfx_0)}=\bar{x}
\end{equation}
Suppose that $\bar{x}>0$. Then it follows from
(\ref{eq:cor:steady_state_2}) that there exists time instant
$t_1<\infty$ and  some constant $0<\delta<\bar{x}$ such that
\[
\norms{\bfx(t)}\geq \delta \ \forall  t\geq t_1
\]
Hence using (\ref{eq:cor:steady_state_1}) and noticing that
$\gamma_1\in\mathcal{K}_e$ we obtain
\[
\infty>h(\bfz_0)-h(\bfz_0)\geq\lim_{T\rightarrow\infty}\int_{t_1}^T\gamma_1(\delta)d\tau=\infty
\]
Thus we obtained a contradiction. Hence,  $\bar{x}=0$ and,
consequently,
\[
\lim_{t\rightarrow\infty}\norms{\bfx(t)}=0
\]
Then, according to the notion of steady-state characteristic in
Definition \ref{defn:steady_state_norm} this is only possible if
$\bar{h}\in\chi^{-1}(0)$. {\it The lemma is proven.}

\vskip 2mm

{\it Proof of Lemma \ref{lem:steady_state_average_convergence}.}
Analogously to the proof of Lemma
\ref{lem:steady_state_convergence} we notice that
(\ref{eq:cor:steady_state_1}) holds. This, however, implies that
for any constant and positive $T$ the following limit
\[
\lim_{t\rightarrow\infty}\int_{t}^{t+T}\gamma_1(\norms{\bfx(\tau)})d\tau
\]
exists and equals zero. Furthermore, $h(\bfz(t,\bfz_0))\in
[0,h(\bfz_0)]$ for all $t\geq t_0$. Hence, there exists a time
instant $t'$ such that
\[
\norms{\bfx(t)}\leq c\cdot h(\bfz_0)+\varepsilon, \ \forall \
t\geq t',
\]
where $\varepsilon>0$ is arbitrary small. Then taking into account
(\ref{eq:converge_steady_state_average_2}) we can conclude that
\begin{equation}\label{eq:cor:steady_state_3}
\lim_{t\rightarrow\infty}\int_{t}^{t+T}\gamma_1(\norms{\bfx(\tau)})d\tau\geq
\bar{\gamma} \int_{t}^{t+T}\norms{\bfx(\tau)}d\tau=0
\end{equation}
Given that (\ref{eq:cor:steady_state_1_1}) holds, system
(\ref{eq:attracting}) has  the steady-state characteristic on
average and that $\chi_T(\cdot)$ has no zeros in the positive
domain, limiting relation (\ref{eq:cor:steady_state_3}) is
possible only if $\bar{h}=0$. Then, according to
(\ref{eq:attracting}),
$\lim_{t\rightarrow\infty}\norms{\bfx(t)}=0$. {\it The lemma is
proven.}

\vskip 2mm

{\it Proof of Corollary
\ref{cor:non_uniform_small_gain_attracting_set}.} As follows from
Theorem \ref{theorem:non_uniform_small_gain}, state
$\bfx(t,\bfx_0)\oplus\bfz(t,\bfz_0)$ converges to the set
$\Omega_a$ specified by (\ref{eq:domain_of_convergence}). Hence
$h(\bfz(t,\bfz_0))$ is bounded.  Then, according to
(\ref{eq:integral}), estimate (\ref{eq:cor:steady_state_1}) holds.
This, in combination with condition
(\ref{eq:integral_practically_ISS}), implies that $\bfz(t,\bfz_0)$
is bounded. In  other words
\[
\bfx(t,\bfx_0)\oplus\bfz(t,\bfz_0)\in\Omega' \ \forall \ t\geq t_0
\]
where $\Omega'$ is a bounded subset in $\Real^n\times\Real^m$.
Applying the Bolzano-Weierstrass theorem we can conclude that for
every point $\bfx_0\oplus\bfz_0\in\Omega_\gamma$ there is an
$\omega$-limit set $\omega(\bfx_0\oplus\bfz_0)\subseteq\Omega'$
(non-empty).

As follows from C3) and Lemma \ref{lem:steady_state_convergence}
the following holds:
\[
\lim_{t\rightarrow\infty}h(\bfz(t,\bfz_0))\in\chi^{-1}(0)
\]
Therefore, given that $h(\cdot)\in\mathcal{C}^0$, we can obtain
that
\[
\lim_{t_i\rightarrow\infty}h(\bfz(t_i,\bfz_0))=h(\lim_{t_i\rightarrow\infty}\bfz(t_i,\bfz_0))=h(\omega_z(\bfx_0\oplus\bfz_0))\in\chi^{-1}(0)
\]
In other words:
\[
\omega_z(\bfx_0\oplus\bfz_0)\subseteq \Omega_h= \{\bfx\in\Real^n,
\ \bfz\in\Real^m| \ h(\bfz)\in\chi^{-1}(0)\}
\]
Moreover
\[
\omega_x(\bfx_0\oplus\bfz_0)\subseteq \Omega_a= \{\bfx\in\Real^n,
\ \bfz\in\Real^m| \ \norms{\bfx}=0\}
\]
According to assumption C1, the flow
$\bfx(t,\bfx_0)\oplus\bfz(t,\bfz_0)$ is generated by a system of
autonomous differential equations with locally Lipschitz
right-hand side. Then, as follows from \cite{Khalil_2002} (Lemma
4.1, page 127)
\[
\lim_{t\rightarrow\infty}\dist(\bfx(t,\bfx_0)\oplus\bfz(t,\bfz_0),\omega(\bfx_0\oplus\bfz_0))=0
\]
Noticing that
\[
\dist(\bfx(t,\bfx_0)\oplus\bfz(t,\bfz_0),\omega(\bfx_0\oplus\bfz_0))\geq
\dist(\bfx(t,\bfx_0),\Omega_a)+\dist(\bfz(t,\bfz_0),\Omega_h)
\]
we can finally obtain that
\[
\lim_{t\rightarrow\infty}\dist(\bfx(t,\bfx_0),\Omega_a)=0, \
\lim_{t\rightarrow\infty}\dist(\bfz(t,\bfz_0),\Omega_h)=0
\]
{\it The corollary is proven.}

\vskip 2mm

{\it Proof of Corollary \ref{cor:non_uniform_small_gain_GAS}.} As
follows from Theorem \ref{theorem:non_uniform_small_gain}, the
corollary will be proven if Conditions \ref{assume:partition_of_z}
-- \ref{assume:phi_system} are satisfied and also
(\ref{eq:tau_system_1}), (\ref{eq:non_uniform_small_gain}), and
(\ref{eq:tau_system_2}) hold. In order to satisfy Condition
\ref{assume:partition_of_z} we select the following sequence
$\mathcal{S}$:
\begin{equation}\label{eq:cor3:sigma_sequence}
\mathcal{S}=\{\sigma_i\}_{i=0}^\infty, \
\sigma_{i}=\frac{1}{\kappa^i}, \ \kappa\in\Real_+, \ \kappa>1
\end{equation}
Let us chose sequences $\mathcal{T}$ and $\Xi$ as follows:
\begin{equation}\label{eq:cor3:tau_sequence}
\mathcal{T}=\{\tau_i\}_{i=0}^\infty, \ \tau_i=\tau^\ast,
\end{equation}
\begin{equation}\label{eq:cor3:xi_sequence}
\Xi=\{\xi_i\}_{i=0}^\infty, \ \xi_i=\xi^\ast,
\end{equation}
where $\tau^{\ast}$, $\xi^\ast$ are positive constants yet to be
defined. Notice that choosing $\mathcal{T}$ as in
(\ref{eq:cor3:tau_sequence}) automatically fulfills condition
(\ref{eq:tau_system_2}) of Theorem
\ref{theorem:non_uniform_small_gain}. On the other hand, taking
into account (\ref{eq:tau_system}),
(\ref{eq:uniformly_attractive_dynamics}) and that $\beta_t(t)$ is
monotonically decreasing in $t$, this choice defines a constant
$\xi^\ast$ as follows:
\begin{equation}\label{eq:cor3:tau_sequence_1}
\beta_t(\tau^\ast)\leq \xi^\ast \beta_t(0) < \beta_t(0), \
0\leq\xi^\ast<1
\end{equation}
Given that the inverse $\beta_t^{-1}$ exists,
(\ref{eq:beta_inverse}), this choice is always possible. In
particular, (\ref{eq:cor3:tau_sequence_1}) will be satisfied for
the following values of $\mathcal{\tau^\ast}$:
\begin{equation}\label{eq:cor3:tau_constant}
\tau^\ast\geq \beta_t^{-1}\left(\xi^\ast\beta_t(0)\right)
\end{equation}
Let us now find the values for $\tau^\ast$ and $\xi^\ast$ such
that Condition \ref{assume:phi_system} is also satisfied. To this
purpose consider systems of functions $\Phi$, $\Upsilon$ specified
by equations  (\ref{eq:phi_system}), (\ref{eq:u_system}). Notice
that function $\beta(s,0)$ in (\ref{eq:phi_system}),
(\ref{eq:u_system}) is linear for system
(\ref{eq:uniformly_attractive_dynamics})
\[
\beta(s,0)=s\cdot\beta_t(0),
\]
and therefore the functions $\rho_{\phi,j}(\cdot)$,
$\rho_{\upsilon,j}$ are identity maps. Hence, $\Phi$, $\Upsilon$
reduce to the following
\begin{equation}\label{eq:cor3:phi_system}
\Phi: \
\begin{array}{ll}
\phi_j(s)&=\phi_{j-1}\cdot \xi^\ast \cdot\beta(s,0)=
\xi^\ast\cdot\beta_t(0)\cdot\phi_{j-1}(s),
\ j=1,\dots,i\\
\phi_0(s)&=\beta_t(0)\cdot s
\end{array}
\end{equation}
\begin{equation}\label{eq:cor3:u_system}
\Upsilon: \
\begin{array}{ll}
\upsilon_j(s)&=\phi_{j-1}(s),
\ j=1,\dots,i\\
\upsilon_0(s)&=\beta_t(0)\cdot s
\end{array}
\end{equation}
Taking into account (\ref{eq:cor3:sigma_sequence}),
(\ref{eq:cor3:phi_system}), (\ref{eq:cor3:u_system}) let us
explicitly formulate requirements (\ref{eq:contraction_1}),
(\ref{eq:contraction_2}) in Condition \ref{assume:phi_system}.
These conditions are equivalent to the boundedness of the
following functions
\begin{equation}\label{eq:cor3:B1}
\norms{\bfx(t_0)}\cdot\beta_t(0)\cdot\kappa^{n} (\xi^\ast\cdot
\beta_t(0))^n;
\end{equation}
\begin{equation}\label{eq:cor3:B2}
\begin{split}
& \kappa^{n}\left(\beta_t(0)\frac{c|h(\bfz_0)|}{\kappa^n}+
\frac{\beta_t(0) c|h(\bfz_0)|}{\kappa^{n-1}}
+\beta_t(0)\sum_{i=2}^n c
|h(\bfz_0)|\frac{1}{k^{n-i}}(\xi^\ast\cdot
\beta_t(0))^{i-1} \right)\\
&=\beta_t(0)c|h(\bfz_0)|+\beta_t(0) c|h(\bfz_0)|\kappa\left(1 +
\sum_{i=2}^n \kappa^{i-1}(\xi^\ast\cdot \beta_t(0))^{i-1}\right)
\end{split}
\end{equation}
Boundedness of the functions  $B_1(\norms{\bfx_0})$ and
$B_2(|h(\bfz_0)|,c)$ is ensured if $\xi^\ast$ satisfies the
following inequality
\begin{equation}\label{eq:cor3:xi_bound}
\xi^\ast\leq \frac{d}{\kappa\cdot\beta_t(0)}
\end{equation}
for some $0\leq d< 1$. Notice that $\kappa>1$, $\beta_t(0)\geq 1$
imply that $\xi^\ast\leq 1$ and therefore constant $\tau^\ast$
satisfying (\ref{eq:cor3:tau_constant}) will always be defied.
Hence, according to (\ref{eq:cor3:B1}), (\ref{eq:cor3:B2}), the
functions $B_1(\norms{\bfx_0})$ and $B_2(|h(\bfz_0)|,c)$
satisfying Condition \ref{assume:phi_system} can be chosen as
\begin{equation}\label{eq:cor3:B1_and_B2}
B_1(\norms{\bfx_0})=\beta_t(0)\norms{\bfx_0}; \
B_2(|h(\bfz_0)|,c)=\beta_t(0)\cdot c\cdot
|h(\bfz_0)|\left(1+\frac{\kappa}{1-d}\right)
\end{equation}

In order to apply Theorem \ref{theorem:non_uniform_small_gain} we
have to check the remaining conditions (\ref{eq:tau_system_1}) and
(\ref{eq:non_uniform_small_gain}). This requires the possibility
of factorization (\ref{eq:assume:gamma_0}) for the function
$\gamma_0(\cdot)$. According to assumption
(\ref{eq:Lipschitz_gamma}) of the corollary the function
$\gamma_0(\cdot)$ is Lipschitz:
\[
|\gamma_0(s)|\leq D_{\gamma,0}\cdot |s|
\]
This allows us to choose function $\gamma_{0,1}(\cdot)$ and
$\gamma_{0,2}(\cdot)$ as follows:
\begin{equation}\label{eq:cor3:gamma_factorization}
\gamma_{0,1}(s)=s, \ \gamma_{0,2}(s)=D_{\gamma,0}\cdot s
\end{equation}
Condition (\ref{eq:tau_system_1}), therefore, is equivalent to
solvability of the following inequality:
\begin{equation}\label{eq:cor3:tau_constant_1}
\left(\frac{1}{\kappa^i}-\frac{1}{\kappa^{i+1}}\right)\frac{\kappa^{i}}{\tau^\ast}\geq
\Delta_0
\end{equation}
Taking into account inequalities (\ref{eq:cor3:tau_constant}),
(\ref{eq:cor3:xi_bound}) we can derive that solvability of
\begin{equation}\label{eq:cor3:Delta_constant}
\Delta_0=
\left(\beta_t^{-1}\left(\frac{d}{\kappa}\right)\right)^{-1}\frac{\kappa-1}{\kappa}
\end{equation}
implies existence of $\Delta_0>0$ satisfying
(\ref{eq:cor3:tau_constant_1}) and, consequently, condition
(\ref{eq:tau_system_1}) of Theorem
\ref{theorem:non_uniform_small_gain}. Given that $d<1$, $\kappa>1$
and $\beta_t(0)\geq 1$ a  positive solution to
(\ref{eq:cor3:Delta_constant}) is always defined. Hence, the proof
will be complete and the claim is non-vacuous if the domain
\begin{equation}\label{eq:cor3:domain}
D_{\gamma,0}\leq
\left(\beta_t^{-1}\left(\frac{d}{\kappa}\right)\right)^{-1}\frac{\kappa-1}{\kappa}
\frac{h(\bfz_0)}{\beta_t(0)\norms{\bfx_0}+\beta_t(0)\cdot c\cdot
|h(\bfz_0)|\left(1+\frac{\kappa}{1-d}\right)+c|h(\bfz_0)|}
\end{equation}
is not empty. {\it The corollary is proven.}

\vskip 2mm

{\it Proof of Corollary \ref{cor:small_gain_like}.} It follows
from Corollary \ref{cor:non_uniform_small_gain_GAS} that state of
the interconnection converges into $\Omega_a$ for all initial
conditions $\bfx_0$, $\bfz_0$ satisfying (\ref{eq:cor3:domain}).
In other words the following inequality should hold:
\begin{equation}\label{eq:cor_small_gain:1}
\begin{split}
&D_{\gamma,0}\left(\beta_t(0)\norms{\bfx_0}+\beta_t(0)\cdot c\cdot
|h(\bfz_0)|\left(1+\frac{\kappa}{1-d}\right)+c|h(\bfz_0)|\right)\leq\\
&\left(\beta_t^{-1}\left(\frac{d}{\kappa}\right)\right)^{-1}
\frac{\kappa-1}{\kappa} \cdot {h(\bfz_0)}
\end{split}
\end{equation}
Hence, assuming that $h(\bfz_0)>0$ we can rewrite
(\ref{eq:cor_small_gain:1}) in the following way:
\begin{equation}\label{eq:cor_small_gain:2}
\begin{split}
D_{\gamma,0}\cdot\beta_t(0)\norms{\bfx_0}&\leq\\
&\left(\left(\beta_t^{-1}\left(\frac{d}{\kappa}\right)\right)^{-1}
\frac{\kappa-1}{\kappa} - D_{\gamma,0}\cdot
c\left(\beta_t(0)\cdot\left(1+\frac{\kappa}{1-d}\right)+1\right)
\right){h(\bfz_0)}
\end{split}
\end{equation}
Solutions to (\ref{eq:cor_small_gain:2}) exist, however, if the
inequality
\[
\left(\beta_t^{-1}\left(\frac{d}{\kappa}\right)\right)^{-1}
\frac{\kappa-1}{\kappa} \geq D_{\gamma,0}\cdot
c\left(\beta_t(0)\cdot\left(1+\frac{\kappa}{1-d}\right)+1\right)
\]
or, equivalently
\begin{equation}\label{eq:cor_small_gain:3}
D_{\gamma,0}\cdot
c\cdot\left(\beta_t(0)\cdot\left(1+\frac{\kappa}{1-d}\right)+1\right)
\cdot \beta_t^{-1}\left(\frac{d}{\kappa}\right)
\frac{\kappa}{\kappa-1}<1
\end{equation}
is satisfied. The estimate of the trapping region follows from
(\ref{eq:cor_small_gain:2}).

Let us finally show that continuity of $h(\bfz)$ implies that the
volume of $\Omega_\gamma$ is nonzero in $\Real^n\oplus\Real^m$.
For the sake of compactness we rewrite inequality
(\ref{eq:cor_small_gain:2}) in the following form:
\begin{equation}\label{eq:cor_small_gain:4}
\norms{\bfx_0}\leq C_{\gamma} h(\bfz_0),
\end{equation}
where $C_\gamma$ is a constant depending on $d$, $\kappa$,
$\beta_t(0)$, and $D_{\gamma,0}$. Given that
(\ref{eq:cor_small_gain:3}) holds we can conclude that
$C_\gamma>0$. According to (\ref{eq:cor_small_gain:4}), domain
$\Omega_\gamma$ contains the following set:
\[
\{\bfx_0\in\Real^n, \ \bfz_0\in\Real^m| \ h(\bfz_0)>
D_z\in\Real_+, \ \norms{\bfx_0}\leq C_\gamma D_z \}
\]

Consider the following domain:
$\Omega_{\bfx,\gamma}=\{\bfx_0\in\Real^n | \ \norms{\bfx_0}\leq
C_\gamma D_z\}$. Clearly, it contains a point
$\bfx_{0,1}\in\Real^n: \ \norms{\bfx_{0,1}}=\frac{C_\gamma
D_z}{2}$. For the point $\bfx_{0,1}$ and for all
$\varevec_1\in\Real^n: \ \|\varevec_1\|\leq \frac{C_\gamma
D_z}{4}$ we have that
$\norms{\bfx_{0,1}+\varevec_1}=\inf_{\bfq\in\mathcal{A}}\|\bfx_{0,1}+\varevec_1-\bfq\|\leq
\inf_{\bfq\in\mathcal{A}}\{\|\bfx_{0,1}-\bfq\|+\|\varevec_1\|\}\leq
\frac{3 C_\gamma D_z }{4}$. On the other hand
$\norms{\bfx_{0,1}+\varevec_1}=\inf_{\bfq\in\mathcal{A}}\|\bfx_{0,1}+\varevec_1-\bfq\|\geq
\inf_{\bfq\in\mathcal{A}}\{\|\bfx_{0,1}-\bfq\| -
\|\varevec_1\|\}\geq \frac{C_\gamma D_z }{4}$.  This implies that
there exists a set of points
$\bfx_{0,2}=\bfx_{0,1}+\varevec_1\in\Real^n$:
$\|\bfx_{0,1}-\bfx_{0,2}\|\leq \frac{C_\gamma D_z}{4}$,
$\bfx_{0,2}\notin\mathcal{A}$, $\norms{\bfx_{0,2}}\leq {C_\gamma
D_z}$.

Consider now the following domain:
$\Omega_{\bfz,\gamma}=\{\bfz_0\in\Real^m| \ h(\bfz_0)> D_z\}$. Let
us pick $\bfz_{0,1}\in\Omega_{\bfz,\gamma}$: $h(\bfz_{0,1})= 2
D_z$. Because $h(\cdot)$ is continuous we have that
\[
\forall \ \varepsilon>0, \ \exists \ \delta>0: \
\|\bfz_{0,1}-\bfz_{0,2}\|<\delta \Rightarrow
|h(\bfz_{0,1})-h(\bfz_{0,2})|< \varepsilon
\]
Let $\varepsilon=D_z$, then $-D_z<h(\bfz_{0,1})-h(\bfz_{0,2})<D_z$
and therefore $h(\bfz_{0,2})>D_z$. Hence there exists a set of
points $\bfz_{0,2}\in\Real^m$: $\|\bfz_{0,1}-\bfz_{0,2}\|<\delta$,
$\bfz_{0,2}\in\Omega_{\bfz,\gamma}$.

Consider the following set
\[
\Omega_{\bfx \bfz,\gamma}=\left\{\bfx'\in\Real^n, \
\bfz'\in\Real^m| \
\|\bfx_{0,1}-\bfx'\|^2+\|\bfz_{0,1}-\bfz'\|^2\leq r^2, \ r=\min
\left\{\delta, \frac{C_\gamma D_z}{4}\right\}\right\}
\]
For all $\bfx_0,\bfz_0 \in \Omega_{\bfx \bfz, \gamma}$ we have
that $\bfx_0\in\Omega_{\bfx,\gamma}$,
$\bfz_0\in\Omega_{\bfz,\gamma}$. Hence, inequality
(\ref{eq:cor_small_gain:4}) holds, and
$\bfx_0\oplus\bfz_0\in\Omega_\gamma$. The volume of the set
$\Omega_{\bfx \bfz,\gamma}$ is defined by the volume of the
interior of a sphere in $\Real^{n+m}$ with nonzero radius. Thus
the volume of $\Omega_\gamma\supset\Omega_{\bfx \bfz, \gamma}$ is
also nonzero. {\it The corollary is proven.}

\vskip 2mm

{\it Proof of Corollary \ref{cor:identifier}.} Let
$\lambdavec(\tau,\lambda_0)$ be a solution of system
(\ref{eq:exo_searching}). Consider it as a function of variable
$\tau$. Let us pick some monotone, strictly increasing function
$\sigma$ such that the following holds
\[
\tau=\sigma(t), \ \sigma:\Real_+\rightarrow\Real_+
\]
Given that  $\etavec(\Omega_\lambda)$ is dense in $\Omega_\theta$,
for any $\thetavec\in\Omega_\theta$ there always exists a vector
$\lambdavec_\theta\in\Omega_\lambda$ such that
$\etavec(\lambdavec_\theta)=\thetavec+\epsilon_\theta$, where
$\|\epsilon_\theta\|$ is arbitrary small. Furthermore,
$\lambdavec(\tau)$ is dense in $\Omega_\lambda$, hence there is a
point $\lambdavec^\ast=\lambdavec(\tau^\ast,\lambdavec_0)$, which
is arbitrarily close to $\lambdavec_\theta$. Consider the
following difference
\[
\bff(\xivec(t),\thetavec)-\bff(\xivec(t),\hat{\thetavec})=\bff(\xivec(t),\thetavec)-f(\xivec(t),\etavec(\lambdavec^\ast))+\bff(\xivec,\etavec(\lambdavec^\ast))-\bff(\xivec,\etavec(\lambdavec(\sigma(t))))
\]
The function $\bff(\cdot)$ is locally bounded and $\etavec(\cdot)$
is Lipschitz, then
\[
\|\bff(\xivec,\thetavec)-\bff(\xivec,\etavec(\lambdavec^\ast))\|\leq
D_f\|\epsilon_\theta\|+\Delta_f=\Delta_\theta+\Delta_f
\]
where $\Delta_\theta$ is arbitrary small. Hence
\begin{equation}\label{eq:difference_nonuniform_1}
\begin{split}
&\|\bff(\xivec,\etavec(\lambdavec^\ast))-\bff(\xivec,\etavec(\lambdavec(\sigma(t))))\|\leq
D_f\|\etavec(\lambdavec^\ast)-\etavec(\lambdavec(\sigma(t)))\|+\Delta_f+\Delta_\theta\\
&\leq D_f \cdot D_\eta
\|\lambdavec^\ast-\lambdavec(\sigma(t))\|+\Delta_f+\Delta_\theta
\end{split}
\end{equation}
Noticing that
$\lambdavec^\ast=\lambdavec(\tau^\ast,\lambdavec_0)=\lambdavec(\sigma(t^\ast),\lambdavec_0)$
and taking into account the Poisson stability of
(\ref{eq:exo_searching}), we can always choose
$\lambdavec^\ast(\sigma^\ast,\lambdavec_0)$ such that
$\sigma^\ast>\sigma(t_0)=\tau_0$ for any $\tau_0\in\Real_+$.
Hence, according to (\ref{eq:difference_nonuniform_1}) the
following estimate holds:
\begin{equation}\label{eq:difference_nonuniform_2}
\begin{split}
&
\|\bff(\xivec,\etavec(\lambdavec^\ast))-\bff(\xivec,\etavec(\lambdavec(\sigma(t))))\|\leq
D_f \cdot D_\eta \|\int_{\sigma(t)}^{\sigma^\ast}
S(\lambdavec(\sigma(\tau)))d\tau
\|+\Delta_f+\Delta_\theta\\
& \leq D_f \cdot D_\eta \cdot \max_{\lambdavec\in\Omega_\lambda}
\|S(\lambdavec)\| |\sigma^\ast-\sigma(t)|=\mathcal{D}\cdot
|\sigma^\ast-\sigma(t)|+\Delta_f+\Delta_\theta, \ \mathcal{D}=D_f
\cdot D_\eta \cdot \max_{\lambdavec\in\Omega_\lambda}
\|S(\lambdavec)\|
\end{split}
\end{equation}
Denoting
$\bfu(t)=\bff(\xivec(t),\thetavec)-\bff(\xivec(t),\hat{\thetavec})+\varevec(t)$
we can now conclude that
\begin{equation}\label{eq:input_decomposition}
\begin{split}
\|\bfu(t)\|&\leq \Delta_\epsilon+\Delta_f+
\|\bff(\xivec(t),\thetavec)-f(\xivec(t),\etavec(\lambdavec^\ast))\|+\mathcal{D}\cdot
|\sigma^\ast-\sigma(t)|\\
&\leq \Delta_\epsilon + 2\Delta_f + \Delta_\theta + D_f
\|\thetavec-\etavec(\lambdavec^\ast)\|+ \mathcal{D}\cdot
|\sigma^\ast-\sigma(t)|
\end{split}
\end{equation}
Notice that due to the denseness of $\lambdavec(t,\lambdavec_0)$
in $\Omega_\lambda$ it is always possible to choose
$\lambdavec^\ast$ such that
\[
D_f \|\thetavec-\etavec(\lambdavec^\ast)\|=D_f
\|\etavec(\lambdavec_\theta)-\etavec(\lambdavec^\ast)\|\leq D_f
D_\eta\|\lambdavec_\theta-\etavec(\lambdavec^\ast)\|\leq
\Delta_\lambda
\]
Hence, according to (\ref{eq:input_decomposition}), we have
\[
\|\bfu(t)\|_{\infty,[t_0,t]}\leq
2\Delta_f+\Delta_\varepsilon+\delta+\mathcal{D}\cdot
\|\sigma^\ast-\sigma(t)\|_{\infty,[t_0,t]}
\]
where the term $\delta>\Delta_\theta+\Delta_\lambda$ can be made
arbitrary small.

Therefore Assumption \ref{assume:observer_error_unperturbed}
implies that the following inequality holds:
\begin{equation}\label{eq:normed_description}
\normss{\bfx(t)}{\Delta(M)}\leq
\beta(t-t_0)\normss{\bfx(t_0)}{\Delta(M)}+c \cdot \mathcal{D}\cdot
\|\sigma^\ast-\sigma(t)\|_{\infty,[t_0,t]}
\end{equation}

Let us now define $\sigma(t)$ as follows
\begin{equation}\label{eq:nonumniform_time_scale}
\sigma(t)=\int_{t_0}^{t}\gamma
\normss{\psi(\bfx(\tau))}{\Delta(M)}d \tau
\end{equation}
Moreover, let us introduce the following notation
\[
h(t)=\sigma^\ast-\sigma(t)= \sigma^\ast-\int_{t_0}^{t}\gamma
\normss{\psi(\bfx(\tau))}{\Delta(M)}d \tau
\]
then for all $t',\ t\geq t_0$, $t\geq t'$ we have that
\[
h(t')-h(t)=\int_{t'}^{t} \gamma
\normss{\psi(\bfx(\tau))}{\Delta(M)}d \tau
\]

Taking into account equation (\ref{eq:difference_nonuniform_1}),
(\ref{eq:difference_nonuniform_2}), equality
\[
\frac{\pd \lambdavec(\sigma(t),\lambdavec_0)}{dt}=\frac{\pd
\sigma(t)}{dt} S(\lambdavec(\sigma(t),\lambdavec_0))=\gamma
\normss{\psi(\bfx(\tau))}{\Delta(M)}
S(\lambdavec(\sigma(t),\lambdavec_0)),
\]
equation (\ref{eq:normed_description}), and denoting
$D_\lambda=c\mathcal{D}$, we can conclude that the following holds
along the trajectories of (\ref{eq:identifier_closed}):
\begin{equation}\label{eq:identifier_closed_norms}
\begin{split}
\normss{\bfx(t)}{\Delta(M)}&\leq
\beta(t-t_0)\normss{\bfx(t_0)}{\Delta(M)}+ D_\lambda
\|h(\tau)\|_{\infty,[t_0,t]}\\
h(t_0)-h(t)&=\int_{t_0}^{t} \gamma
\normss{\psi(\bfx(\tau))}{\Delta(M)}d \tau
\end{split}
\end{equation}
Hence, according to Corollary
{\ref{cor:non_uniform_small_gain_attracting_set}}, the limit
relation (\ref{eq:identifier_convergence}) holds for all
$|h(t_0)|$, $\normss{\bfx(t_0)}{\Delta(M)}$ which belong to the
domain
\[
\Omega_\gamma: \ \gamma\leq
\left(\beta_t^{-1}\left(\frac{d}{\kappa}\right)\right)^{-1}\frac{\kappa-1}{\kappa}
\frac{h(t_0)}{\beta_t(0)\normss{\bfx(t_0)}{\Delta+\delta}+\beta_t(0)\cdot
D_\lambda\cdot
|h(t_0)|\left(1+\frac{\kappa}{1-d}\right)+D_\lambda|h(t_0)|}
\]
for some $d<1$, $\kappa>1$. Notice, however, that
$\normss{\bfx(t)}{\Delta+\delta}$ is always bounded as
$\bff(\cdot)$ is Lipschitz in $\theta$ and both $\thetavec$ and
$\hat{\thetavec}$ are bounded ($\etavec(\cdot)$ is Lipschitz and
$\lambdavec(t,\lambdavec_0)$ is bounded according to assumptions
of the corollary). Moreover, due to the Poisson stability of
(\ref{eq:exo_searching}) it is always possible to choose a point
$\lambdavec^\ast$ such that $h(t_0)=\sigma^\ast$ is arbitrary
large. Hence the choice of  $\gamma$ in
(\ref{eq:identifier_closed_norms}) as (\ref{eq:gamma_identifier})
suffices to ensure that $h(t)$ is bounded. Moreover, it follows
that $h(t)$ converges to a limit as $t\rightarrow\infty$. This
implies that $\gamma\int_{t_0}^t\normss{\bfx(\tau)}{\Delta(M)}$
also converges as $t\rightarrow\infty$, and,  consequently,
$\lambdavec(t,\lambdavec_0)$ converges to some
$\lambdavec'\in\Omega_\lambda$. Hence the following holds
\[
\lim_{t\rightarrow\infty}\hat{\thetavec(t)}=\thetavec'
\]
for some $\thetavec'\in\Omega_\theta$. According to the corollary
conditions, system (\ref{eq:observer_unperturbed}) has steady
state characteristics with respect to $\hat{\thetavec}$. Then, in
the same way as  in the proof  of  Lemma
\ref{lem:steady_state_convergence}, we can show that
(\ref{eq:identifier_convergence}) holds. {The corollary is
proven.}

\bibliographystyle{plain}
\bibliography{non_uniform_small_gain}

\end{document}